\documentclass[11pt,a4paper]{amsart}

\usepackage{makecell}
\usepackage{mathpazo}
\usepackage[mathscr]{eucal}
\usepackage{enumerate}
\usepackage[psamsfonts]{amssymb}
\usepackage{pdfsync}
\usepackage{graphicx,psfrag}
\usepackage[psamsfonts]{amssymb}
\usepackage{amscd}
\usepackage{amsmath}
\usepackage{amsxtra}
\usepackage{mathrsfs}
\usepackage{stmaryrd}
\usepackage{comment}
\usepackage{enumerate}
\usepackage{hyperref}
\usepackage[capitalize,nameinlink]{cleveref}

\usepackage{esint}

\usepackage{subcaption}

\usepackage[width=5.8in, height=8.5in, bottom=1.3in, centering]{geometry}

\usepackage[all]{xy}

\usepackage{mathtools}

\usepackage{mathrsfs}

\usepackage{hyperref}
\hypersetup{ colorlinks, citecolor=green, filecolor=black, linkcolor=blue, urlcolor=blue } 

\theoremstyle{plain}
        \newtheorem{theorem}{Theorem}[section]
        \newtheorem{lemma}[theorem]{Lemma}
        \newtheorem{proposition}[theorem]{Proposition}
        \newtheorem{corollary}[theorem]{Corollary}

        \theoremstyle{definition}
        \newtheorem{definition}[theorem]{Definition}
        \newtheorem{remark}[theorem]{Remark}

\newcommand{\pvint}{\mathop{\mathrlap{\pushpv}}\!\int}
\newcommand{\pushpv}{\mathchoice
  {\mkern5mu\rule[.6ex]{.5em}{1pt}}
  {\mkern2.8mu\rule[.5ex]{.35em}{.8pt}}
  {\mkern2.5mu\rule[.29ex]{.3em}{.7pt}}
  {\mkern2mu\rule[.2ex]{.2em}{.5pt}}
}

\newcommand{\calA}{\mathcal{A}}

\newcommand{\C}{\mathbb{C}}
\newcommand{\Z}{\mathbb{Z}}

\newcommand{\Hom}{\operatorname{Hom}}

\newcommand{\cmdtext}{\mathsf}
\newcommand{\EZ}{\cmdtext{EZ}}
\newcommand{\Sh}{\cmdtext{Sh}}
\newcommand{\diag}{\cmdtext{diag}}
\newcommand{\Tot}{\cmdtext{Tot}}
\newcommand{\EZpert}{\cmdtext{EZ}^{\rm{pert}}}
\newcommand{\Nabla}{\nabla}
\newcommand{\Tr}{\cmdtext{Tr}}
\newcommand{\Ch}{\cmdtext{Ch}}
\newcommand{\HKR}{\cmdtext{HKR}}
\newcommand{\Ad}{\cmdtext{Ad}}
\newcommand{\End}{\cmdtext{End}}
\newcommand{\cs}{\cmdtext{cs}}
\newcommand{\pr}{\cmdtext{pr}}
\newcommand{\id}{\cmdtext{id}}
\newcommand{\dR}{\cmdtext{dR}}
\newcommand{\tr}{\Tr}
\newcommand{\Sym}{\cmdtext{Sym}}
\newcommand{\ev}{\cmdtext{ev}}

\newcommand{\inv}{\cmdtext{inv}}

\renewcommand{\AA}{\mathbf{A}}

\newcommand{\sL}{L}
\newcommand{\sC}{C}
\newcommand{\Hoch}{\cmdtext{Hoch}}
\newcommand{\sH}{H}
\newcommand{\sHP}{HP}
\newcommand{\sCC}{CC}

\renewcommand{\tilde}{\widetilde}

\author{Bjarne Kosmeijer}
\address{University of Amsterdam }
\email{b.a.kosmeijer@uva.nl}
\author{Hessel Posthuma}
\address{University of Amsterdam }
\email{h.b.posthuma@uva.nl}

\title{A note on the equivariant Chern character in Noncommutative Geometry}
\begin{document}
\maketitle
\begin{abstract}
Given a smooth action of a Lie group on a manifold, we give two constructions of the Chern character of an equivariant vector bundle in the 
cyclic cohomology of the crossed product algebra. The first construction associates a cycle to the vector bundle whose structure maps are closely related to Getzler's model for equivariant cohomology. The second construction uses a direct map between this model and the (periodic) cyclic cohomology localized at the unit element. Finally, it is shown that the two constructions are equivalent when the action is proper. 
\end{abstract}

\tableofcontents

\section*{Introduction}
Noncommutative geometry aims to study the geometry of certain singular spaces that are modeled using noncommutative algebras, generalizing ordinary geometry when applied to the commutative case. A prime example of such a noncommutative geometry is the singular space defined by the quotient of an action of a Lie group $G$ on a manifold $M$, which is modeled using the crossed product algebra of the induced action of $G$ on the ring of smooth functions $C^\infty_c(M)$.

Cyclic homology is a homology theory for algebras that plays the role of de Rham cohomology of manifolds in noncommutative geometry. The cyclic homology and 
cohomology of the crossed product algebra $\calA:=G\ltimes C^\infty_c(M)$ and its connection to equivariant cohomology has been extensively studied from the 
early days of cyclic homology on:
\begin{itemize}
\item[(I)] Connes \cite{connes-transverse} and Feigin--Tsygan \cite{FT} studied the case of the action of a discrete group, followed by Getzler--Jones \cite{gj}, Nistor \cite{Nistor-discrete}, and Ponge \cite{Ponge-ga,Ponge-proc}. 
In these papers the cyclic homology of the crossed product algebra is related to the equivariant cohomology of fixed points of the action (``twisted sectors'').
\item[(II)] Brylinski \cite{Brylinski-eq} and Block--Getzler \cite{BG} considered action of compact Lie groups: the cyclic homology is given as a the cohomology of a sheaf over the group whose 
stalk at each element is given in terms of the Cartan model for equivariant cohomology of its fixed point set.
\item[(III)] Finally, Block--Getzler--Jones \cite{BGJ} and Nistor \cite{nistor} generalized the previous case to the setting where the group is allowed to be non-compact. This case is reduced to 
the compact case after localization.
\end{itemize}
These results are closely related to computations for \'etale groupoids (see Brylinsky--Nistor \cite{BN} and Crainic \cite{Crainic}) for (I) and proper Lie groupoids (see \cite{ppt}) for (II).

In the first two cases above the isomorphisms can be explicitly constructed on the chain-level, generalizing the classical HKR map, and can therefore be used to 
compare Connes' noncommutative Chern character of vector bundles in cyclic homology with the classical equivariant Chern character in the equivariant cohomology.
In the third case this is more problematic because the isomorphisms are not explicitly given on the chain-level. 
In \cite{gorokhovsky}, Gorokhovsky studies the Chern character in noncommutative geometry for actions of discrete groups, using the formalism of cycles in the dual cyclic cohomology, and proves compatibility with the equivariant Chern character under the dual isomorphism of (I), localized at the identity, above.
In this paper we generalize this result to the general case of a (not necessarily compact) unimodular Lie group acting on an oriented manifold. 

The results of this paper can be summarized by the diagram:
\[
\xymatrix{{\rm Vect}_G(M)\ar[d]_{{\rm Ch}}\ar[rr]^{{\rm Ch}_{\rm NC}}&&HP^{\bullet}(\mathcal{A})\\H^\bullet_G(M)\ar[urr]_{\Phi}&&}
\]
In this diagram, ${\rm Vect}_G(M)$ denotes the set of isomorphism classes of equivariant vector bundles and the vertical arrow is the classical equivariant Chern character 
in $H^\ev_G(M)=H^{\bullet}(EG\times_GM)$.
In this paper, Getzler's model \cite{getzler} for equivariant cohomology, generalizing Cartan's model for compact Lie groups, plays an important role: 
First, in \cref{sec:ccc} we define the horizontal arrow in the diagram by constructing cycles over the crossed product algebra whose structure maps are closely related to the 
differentials of his complex. Second, we construct the diagonal map $\Phi$ in \cref{sec:cec} by generalizing the HKR maps of \cite{BG} localized at the identity
 to the non-compact case and showing that it pairs naturally with cocycles in Getzler's model. This map generalizes the map in \cite{connes-transverse} for actions of discrete groups.  Finally, in \cref{chernsagree} we prove that the diagram commutes when
the action is proper. For non-proper actions, commutativity of the diagram remains open. Let us also remark that the theory discussed in this paper is localized at the unit elements, the extension to other localizations and related fixed points is left for future research
\subsection*{Acknowledgements} 
The research of B.K. is supported by the Dutch Research Council (NWO) through project nr. 613.001.021.
\section{The Chern character via cycles}
\label{sec:ccc}
\subsection{Preliminaries}
Let $G$ be a Lie group acting smoothly (from the right) on a manifold $M$. We assume that $G$ is unimodular and that $M$ is oriented. 
Associated to this data is the Lie groupoid over $M$ given by $\mathsf{G}:=M\times G$ and source and target maps given by $s(x,g)=xg$ and 
$t(x,g)=x$. The convolution algebra of this groupoid is given by $\mathcal{A}:=C^\infty_c(\mathsf{G})=C^\infty_c(M\times G)$ equipped with the 
product
\begin{equation}
\label{conv}
(f_1*f_2)(x,g):=\int_Gf_1(x,h)f_2(xh,h^{-1}g)dh.
\end{equation}

Let us briefly recall the definition of the cyclic cohomology of $\mathcal{A}$, c.f.\ \cite{Loday}. Denote by $\mathcal{A}^+$ the algebra 
given by adjoining a unit to $\mathcal{A}$, and by $C^k(\mathcal{A}^+)$ the space of $(k+1)$-multilinear forms on $\mathcal{A}^+$ that satisfy
\[
\phi(a_0,\ldots,a_{i-1},1,a_i,\ldots,a_{k-1})=0,\qquad\mbox{for all}~i=1,\ldots,k.
\]
The Hochschild differential maps $b:C^k(\mathcal{A}^+)\to C^{k+1}(\mathcal{A}^+)$ whereas Connes' differential maps 
$B:C^k(\mathcal{A}^+)\to C^{k-1}(\mathcal{A}^+)$ satisfying $b^2=0=B^2$ and $[b,B]=0$. We then set
\[
CC^k(\mathcal{A}^+):=\bigoplus_{n\in\Z} C^{k+2n}(\mathcal{A}^+),
\]
equipped with the differential $b+B:CC^k(\mathcal{A}^+)\to CC^{k+1}(\mathcal{A}^+)$. The cyclic cohomology $HC^\bullet(\mathcal{A})$ of 
$\mathcal{A}$ is then defined as the cohomology of the subcomplex given by the kernel of the canonical morphism 
$CC^\bullet(\mathcal{A}^+)\to CC^\bullet(\C)$.
\begin{remark}
The crossed product algebra $\mathcal{A}:=C^\infty_c(G\times M)$ has a natural locally convex topology and we shall therefore actually 
use (inductive) completed tensor products $\hat{\otimes}$ so that $\mathcal{A}^{\hat{\otimes} (k+1)}\cong C^\infty_c(G^{\times (k+1)}\times M^{\times(k+1)})$, and continuous homomorphisms: $C^k(\mathcal{A})={\rm Hom}_{\rm cont}(\mathcal{A}^{\hat{\otimes}(k+1)},\C)$. With this
definition, the cyclic cohomology of $\mathcal{A}$ can be given the structure of a module over the ring $C^\infty_{\rm inv}(G)$ of 
conjugacy invariant functions given by $(f\phi)(a_0\otimes\ldots\otimes a_k):=\phi(f \cdot (a_0\otimes\ldots\otimes a_k))$ where $f\in C^\infty_{\rm inv}(G)$ acts by
\[
(f \cdot (a_0\otimes\ldots\otimes a_k))(g_0,\ldots,g_k,x_0,\ldots,x_k):=f(g_0\cdots g_k)a_0(g_0,x_0)\cdots a_k(g_k,x_k).
\]
This module structure is compatible with the differentials $b$ and $B$ and therefore gives the cyclic cohomology $HC(\mathcal{A})$ the 
structure of a $C^\infty_{\rm inv}(G)$-module. The dual of this module structure was used in \cite{nistor}  to analyse the cyclic homology 
of $\mathcal{A}$ by localization to the maximal ideals $\mathfrak{m}_{(g)}$ of functions in $C^\infty_{\rm inv}(G)$ vanishing at the 
conjugacy class $(g)$ of elements $g\in G$. In this paper we shall be mostly concerned with localization at the unit element $e\in G$.
\end{remark}

In this section we shall construct cyclic cohomology classes for the convolution algebra using Connes' formalism of \textit{cycles}, c.f. \cite[\S 3.1]{connes-book}. The ingredients for this construction are an (externally) curved DGA, together with a closed graded trace. We now briefly recall these concepts. 
\begin{definition}\label{gcdga}
An \textit{externally curved DGA} is a graded vector space $\Omega=\Omega^\bullet$ with an associative graded product $\ast$, a differential $D:\Omega^\bullet\to\Omega^{\bullet+1}$ and a multiplier $\Theta$ on $\Omega$ of degree $2$, such that
\begin{itemize}
\item $D$ is a graded derivation,
\item $D^2=[\Theta,-]$,
\item $D(\Theta\ast\alpha)=\Theta\ast(D\alpha)$.
\end{itemize}
\end{definition}
\begin{remark}{\ }
\begin{itemize}
\item[$i)$]
Recall that a multiplier $\Theta$ of degree $d$ is a pair of linear maps $\Theta_l,\Theta_r:\Omega^\bullet\to\Omega^{\bullet+d}$ satisfying
\begin{itemize}
\item $\Theta_l(\alpha\ast\beta)=\Theta_l(\alpha)\ast\beta$
\item $\Theta_r(\alpha\ast\beta)=\alpha\ast\Theta_r(\beta)$
\item $\Theta_r(\alpha)\ast\beta=\alpha\ast\Theta_l(\beta)$
\end{itemize}
When $\Omega$ is unital, multipliers of degree $d$ are one-to-one with elements $\Theta\in \Omega^d$ where $\Theta_l=\Theta\ast-$ and $\Theta_r=-\ast\Theta$. Due to this $\Theta_l$ and $\Theta_r$ are also in the non-unital case written as $\Theta\ast -$ and $-\ast\Theta$.
\item[$ii)$]
We will induce multipliers of degree $d$ on a graded algebra $\Omega$ by injecting $\Omega$ into a bigger graded algebra and find an element $\theta$ of degree $d$ such that $\theta\Omega\subset\Omega$ and $\Omega\theta\subset\Omega$ hold in the bigger algebra. 
\item[$iii)$]
The last property in Definition \ref{gcdga} is called the (left) Bianchi-identity, and corresponds to the fact that $D\Theta=0$ in the unital case,
because in the non-unital case $D\Theta$ may not be defined. There is also a right Bianchi-identity $D(\alpha\ast\Theta)=(D\alpha)\ast\Theta$, but under the assumption that $D^2\alpha=[\Theta,\alpha]$ the left and right Bianchi-identities are equivalent since
\begin{equation*}
D(\Theta\ast\alpha-\alpha\ast\Theta)=D(D^2\alpha)=D^2(D\alpha)=\Theta\ast(D\alpha)-(D\alpha)\ast\Theta.
\end{equation*}
\end{itemize}
\end{remark}
\begin{definition}
Let $(\Omega,\ast,D,\Theta)$ be an externally curved DGA. A \textit{closed graded trace of degree $n$} on $\Omega$ is a functional $\pvint:\Omega^n\to\mathbb{C}$ such that
\begin{itemize}
\item $\pvint \alpha\ast\beta=(-1)^{|\alpha||\beta|}\pvint \beta\ast\alpha$,
\item $\pvint \Theta\ast\alpha=\pvint \alpha\ast\Theta$,
\item $\pvint D\alpha=0$.
\end{itemize}
\end{definition}
Suppose now given a triple $(\Omega,\pvint,\rho)$, where $\Omega$ is an externally curved DGA, $\pvint$ a closed graded trace of degree $n$ and $\rho:\mathcal{A}\to\Omega^0$ a morphism of associative algebras.
In \cite[\S 3.1.1]{connes-book}, Connes constructs a cyclic cocycle on $\mathcal{A}$ out of such a triple.
Here we recall Gorokhovsky's JLO-type formula \cite[\S 2]{gorokhovsky} for this class. Write $\Delta^k$ for the standard simplex
\begin{equation*}
\Delta^k=\{t_0,...,t_k\geq 0: t_0+\cdots+t_k=1\}
\end{equation*}
with measure $dt_1\cdots dt_k$.
\begin{theorem}\cite[Thm. 2.1.]{gorokhovsky}
\label{thmg}
Let $(\Omega,\ast,D,\Theta)$ be a externally curved DGA with a closed graded trace $\pvint$ of degree $n$, and let $\rho:\mathcal{A}\to\Omega^0$ be a map of associative algebras. The maps ${\rm Ch}_{\Omega,\rho}^k:\mathcal{A}^{\otimes (k+1)}\to\mathbb{C}$ for $0\leq k\leq n$ with $k\equiv n$ mod $2$ defined by
\begin{equation*}
{\rm Ch}^k_{\Omega,\rho}(a_0,...,a_k)=\int_{\Delta^k}\pvint\rho(a_0)\ast e^{-t_0\Theta}\ast D(\rho(a_1))\ast e^{-t_1\Theta}\ast\cdots\ast D(\rho(a_k))\ast e^{-t_k\Theta}dt_1\cdots dt_k
\end{equation*}
satisfy $b{\rm Ch}^k_{\Omega,\rho}=B{\rm Ch}^{k+2}_{\Omega,\rho}$ and hence define a cocycle ${\rm Ch}_{\Omega,\rho}$ of the $(b,B)$-bicomplex $CC^\bullet(\mathcal{A})$.
\end{theorem}
The resulting cyclic cohomolology class in $HC^\bullet(\mathcal{A})$ is called the {\em Chern character} of the cycle $(\Omega,\pvint,\rho)$.

\subsection{The fundamental cycle}
We shall now construct a externally curved DGA $\Omega$ playing the role of fundamental cycle over the convolution algebra $\mathcal{A}=G\ltimes C^\infty_c(M)$. The elements of $\Omega$ are given by
\[
\Omega:=C^\infty_c\left(G,{\rm Sym}(\mathfrak{g}^*)\otimes\Omega_c(M)\right)
\]
In the following we shall write elements $\alpha\in\Omega$ as maps
\[
(g,X)\mapsto \alpha(g,X)\in\Omega_c(M),
\]
that are smooth and compactly supported in $g\in G$ and polynomial in $X\in\mathfrak{g}$.

On $\Omega$ we introduce $3$ structures: the differential $D$, multiplication $\ast$ and curvature $\Theta$:
\begin{itemize}
\item The differential $D:=d_{\rm dR}+\iota$ is the differential of equivariant cohomology acting on ${\rm Sym}(\mathfrak{g}^*)\otimes\Omega_c(M)$, with $d_{\rm dR}$ given by
\begin{equation*}
(d_{\rm dR}\alpha)(g,X)=d_{\rm dR}(\alpha(g,X))
\end{equation*}
and $\iota$ given by
\begin{equation*}
(\iota\alpha)(g,X)=\iota_{X_M}(\alpha(g,X))
\end{equation*}
\item The multiplication is given by the following generalization of the convolution product \eqref{conv}:
\begin{equation}
\label{conv-df}
(\alpha *\beta)(g,X):=\int_G\alpha(h,X)\wedge h^*\beta(h^{-1}g,{\rm Ad}_{h^{-1}}(X))dh.
\end{equation}
\item The curvature is given by the distribution $\Theta\in C^{-\infty}_c\left(G,\End({\rm Sym}(\mathfrak{g}^*)\otimes\Omega_c(M))\right)$ defined by
\[
\left<\Theta,\alpha\right>(X):=\left.\frac{d}{dt}\right|_{t=0}\alpha(e^{tX},X).
\]
\end{itemize}
The degree of an element $\alpha\in\Omega$ is given, as in equivariant cohomology, by adding the degree of $\alpha$ as a differential form and twice the degree of the polynomial on $\mathfrak{g}$. With this,
the differential $D$ indeed has degree $1$ and the curvature $\Theta$ is of degree $2$.

Note that the convolution product \eqref{conv-df} is essentially the convolution product induced by the $G$-algebra ${\rm Sym}(\mathfrak{g}^\ast)\otimes\Omega_c(M)$ with product given by the pointwise wedge product and the $G$-action given by:
\begin{equation*}
(g\cdot \alpha)(X)=g^\ast(\alpha({\rm Ad}_{g^{-1}}(X)))
\end{equation*}
\begin{remark}\label{conv-distr}
By a distribution 
\[
T\in C^{-\infty}_c\left(G,\End({\rm Sym}(\mathfrak{g}^*)\otimes\Omega_c(M))\right)
\]
we mean a $t$-fibered distribution on the Lie groupoid $G\times\rightrightarrows M$ as in \cite{lmv}, i.e., a $C^\infty(M)$-linear map
\[
T:C^\infty_c(G,{\rm Sym}(\mathfrak{g}^*)\otimes\Omega_c(M))\to {\rm Sym}(\mathfrak{g}^*)\otimes\Omega_c(M),
\]
that is continuous for the Fr\'ech\`et topology. The extension of the convolution product \eqref{conv-df} to distributions is then defined by
\[
\left<T_1* T_2,\varphi\right>:=\left<T_1(g_1),g_1\cdot \left<T_2(g_2),g_1^{-1}\cdot\varphi(g_1g_2)\right>\right>,\qquad\varphi\in C^\infty_c(G,{\rm Sym}(\mathfrak{g}^*)\otimes\Omega_c(M)).
\]
Here $g_1,g_2\in G$ are placeholder variables and the equation should be read as follows: fixing $g_1$, 
the map $g_2\mapsto g_1^{-1}\cdot\phi(g_1g_2)$ lies in $C^\infty_c(G,{\rm Sym}(\mathfrak{g}^*)\otimes\Omega_c(M))$ 
and hence can be paired with $T_2$, leading to an element in ${\rm Sym}(\mathfrak{g}^*)\otimes\Omega_c(M)$. 
This yields a map $g_1\mapsto g_1\cdot\langle T_2(g_2),g_1^{-1}\cdot\phi(g_1g_2)\rangle$ which is in 
$C^\infty_c(G,{\rm Sym}(\mathfrak{g}^*)\otimes\Omega_c(M))$ which can then in turn be paired with $T_1$.

Indeed, for distributions $T_\alpha$, defined by $\alpha\in\Omega$, of the form 
\[
\left<T_\alpha,\varphi\right>(X):=\int_G\alpha(g,X)\wedge \varphi(g,X)dg,\qquad\varphi\in C^\infty_c(G,{\rm Sym}(\mathfrak{g}^*)\otimes\Omega_c(M)),
\]
convolution of distributions leads exactly to the formula \eqref{conv-df}.
\end{remark}

\begin{proposition}
$(\Omega,m,D,\Theta)$ is an externally curved DGA in the sense of Definition \ref{gcdga}.
\end{proposition}
\begin{proof}
It is easy to check that the product $m$ is associative.
Let us then show that $D$ is a derivation for this product. First:
\begin{align*}
\iota m(\alpha,\beta)(g,X)&=\int_G\iota_{X_M}\alpha(h,X)\wedge h^*(\beta(h^{-1}g,{\rm Ad}_{h^{-1}}(X)))dh\\
&\hspace{2cm}+(-1)^{|\alpha|}\int_G\alpha(h,X)\wedge \iota_{X_M}h^*(\beta(h^{-1}g,{\rm Ad}_{h^{-1}}(X)))dh\\
&=m(\iota\alpha,\beta)+(-1)^{|\alpha|}\int_G\alpha(h,X)\wedge h^*\iota_{{\rm Ad}_{h^{-1}}(X)_M}(\beta(h^{-1}g,{\rm Ad}_{h^{-1}}(X)))dh\\
&=m(\iota\alpha,\beta)+(-1)^{|\alpha|}m(\alpha,\iota\beta).
\end{align*}
For the other part of the differential $D$, the derivation property follows immediately from the standard properties of the de Rham differential and the wedge product of forms.

By Cartan's magic formula we have $D^2=\mathcal{L}$ with
\[
{\mathcal L}(\alpha)(g,X):=\mathcal{L}_{X_M}\alpha(g,X),
\]
and this equals $[\Theta,-]$ because
\[
(\Theta*\alpha)(g,X)=\left.\frac{d}{dt}\right|_{t=0} e^{tX}\cdot\alpha(e^{-tX}g,X)
\qquad
\mbox{and} 
\qquad
(\alpha*\Theta)(g,X)=\left.\frac{d}{dt}\right|_{t=0}\alpha(e^{-tX}g,X).
\]

The Bianchi identities $D(\Theta\ast\alpha)=\Theta\ast D(\alpha)$, and $D(\alpha\ast\Theta)=D(\alpha)\ast\Theta$ are easily checked 
to hold true, for the former with the observation that $\mathcal{L}_{X_M}$ and $\iota_{X_M}$ commute as operators on $\Omega_c(M)$ since $\mathcal{L}_{X_M}(X_M)=0$.
This proves the proposition.
\end{proof}

Let us now introduce the functional $\pvint:\Omega^n\to\C$, where $n=\dim M$, by the formula
\[
\pvint\alpha:=\int_M\alpha(e,0).
\]
\begin{proposition}\label{simplesttrace}
The functional $\pvint$ is a {\em closed graded trace} on $\Omega$.
\end{proposition}
\begin{proof}
It is straightforward to check that $\pvint$ is a graded trace:
\[
\pvint(\alpha*\beta)=(-1)^{|\alpha| |\beta|} \pvint(\beta*\alpha).
\]
We outline the steps anyway, to emphasis the differences with the graded traces that will come up later in this paper. Starting from the left, we have
\begin{equation*}
\pvint(\alpha\ast\beta)=\int_M\int_G\alpha(g,0)\wedge g^\ast\beta(g^{-1},0)dg.
\end{equation*}
Using the commutation relations in $\Omega(M)$ we see that
\begin{equation*}
\int_M\int_G\alpha(g,0)\wedge g^\ast\beta(g^{-1},0)dg=(-1)^{|\alpha||\beta|}\int_M\int_G g^\ast\beta(g^{-1},0)\wedge\alpha(g,0)dg.
\end{equation*}
Using that integration over $M$ is $G$-invariant we get
\begin{equation*}
(-1)^{|\alpha||\beta|}\int_M\int_G g^\ast\beta(g^{-1},0)\wedge\alpha(g,0)dg=(-1)^{|\alpha||\beta|}\int_M\int_G \beta(g^{-1},0)\wedge (g^{-1})^\ast\alpha(g,0)dg,
\end{equation*}
and lastly, since $G$ is unimodular we may replace the integral over $g$ with one over $g^{-1}$ to obtain
\begin{align*}
(-1)^{|\alpha||\beta|}\int_M\int_G \beta(g^{-1},0)\wedge (g^{-1})^\ast\alpha(g,0)dg&=(-1)^{|\alpha||\beta|}\int_M\int_G\beta(g,0)\wedge g^\ast\alpha(g^{-1},0)dg\\
&=(-1)^{|\alpha||\beta|}\pvint\beta\ast\alpha
\end{align*}
Likewise, the fact that it is closed follows immediately from de Rham's theorem:
\[
\pvint(D\alpha)=\int_Md\alpha(e,0)=0.
\]
The identity $\pvint\Theta\ast\alpha=\pvint\alpha\ast\Theta$ follows trivially from the fact that $(\Theta\ast\alpha)(e,0)=(\alpha\ast\Theta)(e,0)=0$. This finishes the proof.
\end{proof}
By Theorem \ref{thmg}, the triple $(\mathcal{A},\Omega,\pvint)$ gives rise to a cyclic cocycle ${\rm Ch}_\Omega$ of degree $n=\dim M$ which is given in the $(b,B)$-complex by the components
\begin{align*}
{\rm Ch}_\Omega^k(a_0,\ldots,a_k):&=\int_{\Delta^k}\pvint a_0*e^{-t_0\Theta}*Da_1*e^{-t_1\Theta}*\cdots *Da_k*e^{-t_k\Theta}dt_1\ldots dt_k.
\end{align*}
where $k=n~{\rm mod}~2$.

Using the fact that $(\alpha\ast\Theta)(g,0)=0$ we see that this cocycle only has contributions for $k=n$, where it can be explicitely written as
\begin{multline*}
{\rm Ch}_\Omega^n(a_0,\ldots,a_n)=\frac{1}{n!}\int_M\int_{G^{\times n}} a_0(h_1)h_1^\ast da_1(h_2)\wedge\cdots\\
\cdots\wedge (h_1\cdots h_{k-1})^\ast da_{k-1}(h_k)\wedge (h_1\cdots h_k)^\ast da_k((h_1\cdots h_k)^{-1})dh_1\cdots dh_k
\end{multline*}
Using that $G$ is unimodular, then we can also write this as
\begin{multline*}
{\rm Ch}_\Omega^n(a_0,\ldots,a_n)=\frac{1}{n!}\int_M\int_{G^{\times n}} a_0((h_1\cdots h_k)^{-1})((h_1\cdots h_k)^{-1})^\ast da_1(h_1)\wedge\cdots\\
\cdots\wedge ((h_{k-1}h_k)^{-1})^\ast da_{k-1}(h_{k-1})\wedge (h_k^{-1})^\ast da_k(h_k)dh_1\cdots dh_k
\end{multline*}

\begin{remark}
Inspired by equivariant cohomology one would be tempted to define a ${\rm Sym}(\mathfrak{g^*})^G$-valued functional by
\[
\alpha\mapsto \int_M\alpha(e,-),
\]
but this fails to be a trace for the convolution product \eqref{conv-df}. The problem is the adjoint action of $G$ on $X\in\mathfrak{g}$ in formula
\eqref{conv-df} for the product, and this explains why we put $X=0$ in the definition above. To capture the higher degree polynomial 
terms of $\alpha\in\Omega$, one can twist the trace by an element $\gamma\in{\rm Sym}(\mathfrak{g})^G$, viewed as an invariant differential 
operator $\mathsf{D}_\gamma$ on the Lie algebra $\mathfrak{g}$. 

It will follow that if one defines
\[
\pvint_\gamma\alpha:=\int_M \mathsf{D}_\gamma(\alpha)(e,0),
\] 
this also defines a closed graded trace on $\Omega$. Remark that in combination with evaluation at $0\in\mathfrak{g}$, the invariants
${\rm Sym}(\mathfrak{g})^G$ identify as the algebra of distributions supported at $0$ in the form of derivatives of the $\delta$-distribution via $\gamma\mapsto \mathsf{D}_\gamma(\delta_0)$.
\end{remark}
\begin{proposition}\label{prop-invariantpolygradedtrace}
For $\gamma\in ({\rm Sym}^q\mathfrak{g})^G$, the functional $\pvint_\gamma$ defines a closed graded trace on $\Omega$ of degree ${\rm dim}(M)+2q$.
\end{proposition}
\begin{proof}
The degree of $\pvint_\gamma$ follows from the fact that to obtain a top-form on $M$ after applying $D_\gamma$ and applying $0\in\mathfrak{g}$ to $D_\gamma(\alpha)$ we need $\alpha$ to be of degree ${\rm dim}(M)$ in the differential form part and of polynomial degree $q$, i.e. we need $\alpha$ to be of degree ${\rm dim}(M)+2q$.

To see that $\pvint_\gamma$ vanishes on graded commutators, we compare $\pvint_\gamma \alpha\ast\beta$ and $(-1)^{|\alpha||\beta|}\pvint_\gamma\beta\ast\alpha$ for the case $\gamma=v_1\odot\cdots\odot v_q\in({\rm Sym}\mathfrak{g})^G$ by an explicit calculation. The expression $\pvint_\gamma\alpha\ast\beta$ will now look like
\begin{equation*}
\pvint_\gamma\alpha\ast\beta=\left.\frac{d}{dt_1}\right|_{t_1=0}\cdots\left.\frac{d}{dt_q}\right|_{t_q=0}\int_M\int_G\alpha(g,\sum_{i=1}^q t_iv_i)\wedge\beta(g^{-1},g^\ast{\rm Ad}_{g^{-1}}(\sum_{i=1}^q t_iv_i))dg
\end{equation*}
After applying the same manipulations as in the proof of Proposition \ref{simplesttrace}, this will equal
\begin{equation*}
\pvint_\gamma\alpha\ast\beta=(-1)^{|\alpha||\beta|}\left.\frac{d}{dt_1}\right|_{t_1=0}\cdots\left.\frac{d}{dt_q}\right|_{t_q=0}\int_M\int_G\beta(g,\sum_{i=1}^q t_i{\rm Ad}_g(v_i))\wedge g^\ast\alpha(g^{-1},\sum_{i=1}^q t_iv_i)dg
\end{equation*}
Now since $v_1\odot\cdots\odot v_q$ is $G$-invariant we may replace $\{v_i\}_{i=1,...,q}$ by $\{{\rm Ad}_{g^{-1}}(v_i)\}_{i=1,...,q}$ at no cost, to see that
\begin{equation*}
\pvint_\gamma\alpha\ast\beta=(-1)^{|\alpha||\beta|}\left.\frac{d}{dt_1}\right|_{t_1=0}\cdots\left.\frac{d}{dt_q}\right|_{t_q=0}\int_M\int_G\beta(g,\sum_{i=1}^q t_iv_i)\wedge g^\ast\alpha(g^{-1},\sum_{i=1}^q t_i{\rm Ad}_{g^{-1}}(v_i))dg
\end{equation*}
and this precisely equals $(-1)^{|\alpha||\beta|}\pvint_\gamma\beta\ast\alpha$.

The argument that $\pvint_\gamma$ is closed is the same as the argument for $\pvint$ before, since the $d$-part of $D$ does not contribute to $\pvint_\gamma\circ D$ by deRham's Theorem, while the $\iota$-part of $D$ does not contribute since the only top-form of the form $\iota\omega$ is the $0$-form.
\end{proof}
\begin{remark}
When writing out the JLO-cocycle of Theorem \ref{thmg} for this closed graded trace, one recognizes, just like in the case that $\gamma=0$, that it only has contributions in degree ${\rm dim}(M)$. Indeed, looking at the contribution for a given $k$, we need to consider $(f_0,...,f_k)\subset C^\infty_c(G\times M)$ and make arguments about
\begin{equation*}
f_0\ast e^{-t_0\Theta}\ast Df_1\ast e^{-t_1\Theta}\ast\cdots\ast Df_k\ast e^{-t_k\Theta}
\end{equation*}
However, one notices that this function will eat $g\in G$ and $X\in\mathfrak{g}$ and spit out a $k$-form on $M$, since $(Df_1)(g,X)=d(f_1(g))$ and applying $\Theta$ does not change the degree of the differential form.

So, no matter the application of $D_\gamma$ or evaluating at $g=e$ and $X=0$, we see that integrating over $M$ only gives a non-trivial contribution when $k={\rm dim}(M)$.

Therefore, the resulting cocycle lives in the image of $S^q: HC^{{\rm dim}(M)}(G\ltimes C^\infty_c(M))\to HC^{{\rm dim}(M)+2q}(G\ltimes C^\infty_c(M))$.
\end{remark}
\subsection{Twisting by an equivariant vector bundle}\label{crvDGA}
For an equivariant vector bundle $E\to M$ with a (not necessarily $G$-invariant) connection $\nabla$, there is a variant of the construction of the previous section. In the case that the group is discrete, this construction is due to Gorokhovsky \cite[3]{gorokhovsky}, and we generalize it to  here to the case of an uni-modular Lie group. In this case $\Omega_E:=C^\infty_c(G,{\rm Sym}(\mathfrak{g}^*)\otimes \Omega_c(M,{\rm End}(E)))$ 
and we change the differential to 
\[
(D_\nabla\alpha)(g,X):=d_{\nabla^\End}(\alpha(g,X))+(-1)^{|\alpha|}\alpha(g,X)\wedge\delta(g)+\iota_{X_M}\alpha(g,X),
\] 
where $\nabla^{\rm End}$ is the induced connection on $\End(E)$ and 
\[
\delta(g):=\nabla-g^*\nabla\in\Omega^1(M,{\rm End}(E)).
\]
The curvature is now given by
\[
\Theta_\nabla:=\Theta+(F(\nabla)+\mu)\delta_e,
\]
where $F(\nabla)\in\Omega^2(M,{\rm End}(E))$ is the ordinary curvature of the connection $\nabla$ and $\mu\in\mathfrak{g}^\ast\otimes{\rm End}(E)$ is the moment of $\nabla$ (c.f. \cite[p.518]{BG}) given by
\begin{equation*}
\mu(X)=\nabla_{X_M}-\mathcal{L}_X
\end{equation*}
The multiplication is given by the same formula as \eqref{conv-df}.
\begin{proposition}\label{OmegaEcDGA}
The quadruple $\Omega_{E,\nabla}=(\Omega_E,\ast,D_\nabla,\Theta_\nabla)$ is an externally curved DGA. 
\end{proposition}
We delay the proof of this proposition to \cref{sec:proof}.
We now consider the functional
\[
\pvint\alpha:=\int_M{\rm tr}_E\alpha(e,0)
\]
where ${\rm tr}_E:\Omega^\bullet(M,{\rm End}(E))\to\Omega^\bullet(M)$ is the application of the matrix trace.
\begin{proposition}
The functional $\pvint$ is a closed graded trace on $\Omega_{E,\nabla}$.
\end{proposition}
\begin{proof}
This is similar to the untwisted case. Doing the same steps as in the untwisted	 case, we obtain:
\begin{equation*}
\pvint (\alpha\ast\beta)=(-1)^{|\alpha||\beta|}\pvint (\beta\ast\alpha)+\int_M\int_G{\rm tr}_E([\alpha(g,0),g^\ast\beta(g^{-1},0)])dg,
\end{equation*}
where we note that ${\rm tr}_E([\alpha(g,0),g^\ast\beta(g^{-1},0)])=0$. To see that it is closed we have
\begin{equation*}
\pvint D_\nabla\alpha=\int_M{\rm tr}_E(d_{\nabla^{\rm End}}(\alpha(e,0)))=\int_M d({\rm tr}_E(\alpha(e,0)))=0.
\end{equation*}
To check that $\pvint \Theta_\nabla\ast\alpha=\pvint\alpha\ast\Theta_\nabla$ we first note that
\begin{align*}
(\Theta_\nabla\ast\alpha)(e,0)&=F(\nabla)\wedge \alpha(e,0)\\
(\alpha\ast\Theta_\nabla)(e,0)&=\alpha(e,0)\wedge F(\nabla),
\end{align*}
and then using the fact that taking the trace over ${\rm End}(E)$ is cyclically invariant we obtain
\begin{equation*}
\pvint(\Theta_\nabla\ast\alpha)=\int_M{\rm tr}_E(F(\nabla)\wedge\alpha(e,0))=\int_M{\rm tr}_E(\alpha(e,0)\wedge F(\nabla))=\pvint(\alpha\ast\Theta_\nabla)
\end{equation*}
and hence $\pvint$ is a closed graded trace.
\end{proof}

The cyclic cocycle that we get from Theorem \ref{thmg} is now given by 
\begin{align}
\label{curvedChern}
{\rm Ch}_{\Omega_{E,\nabla}}^k(a_0,\ldots,a_k):=\int_{\Delta^k}\pvint {\rm tr}_E\left(a_0*e^{-t_0\Theta_\nabla}*D_\nabla a_1*e^{-t_1\Theta_\nabla}*\cdots\right.&
\\
\nonumber
\left.\cdots  *D_\nabla a_k*e^{-t_k\Theta}\right)&dt_1\ldots dt_k,
\end{align}
We can also write this out explicitly using the fact that \[(\alpha\ast\Theta_\nabla)(g,0)=\alpha(g,0)\wedge g^\ast F(\nabla),\]
\[(D_\nabla\alpha)(g,0)=d_{\nabla^{\rm End}}(\alpha(g,0))+(-1)^{|\alpha|}\alpha(g,0)\wedge\delta(g),\] and the fact that for $a\in C^\infty_c(M)\subset\Omega^0_c(M,{\rm End}(E))$ we have $d_{\nabla^{\rm End}}a=da$. We end up with
\begin{align*}
{\rm Ch}_{\Omega_{E,\nabla}}^k(a_0,\ldots,a_k)=&\sum_{i_0+\cdots+i_k=\frac{n-k}{2}}\int_{\Delta^k}\frac{(-1)^{\frac{n-k}{2}}t_0^{i_0}\cdots t_k^{i_k}}{i_0!\cdots i_k!}dt
\\
&\int_{G^k}\int_M{\rm tr}_E(a_0(h_1)(h_1^\ast F(\nabla))^{\wedge i_0}\wedge h_1^\ast(da_1(h_2)+a_1(h_2)\delta(h_2))\wedge\cdots
\\
\cdots\wedge(h_1\cdots h_k)^\ast&(da_k((h_1\cdots h_k)^{-1})+a_k((h_1\cdots h_k)^{-1})\delta((h_1\cdots h_k)^{-1})\wedge F(\nabla)^{\wedge i_k})dh,
\end{align*}
where $dt=dt_0\cdots dt_k$ and $dh=dh_1\cdots dh_k$.
We can compute the integral over $\Delta^k$ (see also \cite{gorokhovsky}) to obtain
\begin{multline}\label{classoftwistedCuDGA}
{\rm Ch}_{\Omega_{E,\nabla}}^k(a_0,\ldots,a_k)=\sum_{i_0+\cdots+i_k=\frac{n-k}{2}}\frac{(-1)^{\frac{n-k}{2}}}{\left(\frac{n+k}{2}\right)!}\int_{G^k}\int_M{\rm tr}_E(a_0(h_1)(h_1^\ast F(\nabla))^{\wedge i_0}\wedge\\
\wedge h_1^\ast(da_1(h_2)+a_1(h_2)\delta(h_2))\wedge\cdots\wedge(h_1\cdots h_k)^\ast(da_k((h_1\cdots h_k)^{-1})+\\
+a_k((h_1\cdots h_k)^{-1})\delta((h_1\cdots h_k)^{-1}))\wedge F(\nabla)^{\wedge i_k})dh
\end{multline}
\begin{remark}
Combining two variations, one can also define a $\pvint_\gamma$ on $\Omega_E$ for some $\gamma\in{\rm Sym}(\mathfrak{g})^G$ analogous to the untwisted case to obtain another closed graded trace.
\end{remark}
We finally remark that the cyclic cohomology class does not depend on the choice of connection:
\begin{theorem}
\label{ncchern}
The cyclic cohomology class ${\rm Ch}_{\Omega_{E,\nabla}}\in HC^\bullet(\mathcal{A})$ does not depend on the choice of connection 
$\nabla$ and defines a map
\[
{\rm Ch}_\Omega:{\rm Vect}_G(M)\to HC^\bullet(\mathcal{A}),
\]
where ${\rm Vect}_G(M)$ denotes the set of isomorphism classes of $G$-vector bundles over $M$.
\end{theorem}
\begin{proof}
The invariance under the choice of the connection follows similar to \cite[Lem 2.7]{gorokhovsky}, by remarking that for two connections $\nabla_0,\nabla_1$, the family $\nabla_t=\nabla_0+t(\nabla_1-\nabla_0)$ induces a cobordant family between $\Omega_{E,\nabla_0}$ and $\Omega_{E,\nabla_1}$.
\end{proof}

\section{Compatibility with the equivariant Chern character}
\label{sec:cec}
Recall the basic setting of a Lie group $G$ acting smoothly on a manifold $M$.
The map described in Theorem \ref{ncchern} is called the {\em noncommutative Chern character}, the adjective indicating that it takes values 
in the cyclic cohomology of the noncommutative algebra $\mathcal{A}$ given by the crossed product of the induced action of $G$ on 
$C^\infty(M)$. The {\em classical} ( or: {\em commutative}) {\em equivariant Chern character} takes values in the equivariant cohomology 
$H^\bullet_G(M)$ and is given by taking the ordinary Chern characterof the induced vector bundle over the homotopy quotient 
$(EG\times M)\slash G$.

The aim of this section is to compare these two Chern characters. For this we first construct a map 
\[
\Phi: H^\bullet_G(M)\to HC^\bullet(\mathcal{A}),
\]
using Getzler's model \cite{getzler} for the cochain complex computing the equivariant cohomology $H^\bullet_G(M)$, together 
with computations in {\em cyclic homology} generalizing the case of actions of {\em compact } Lie groups considered in \cite{brylinski,BG,ppt}. We then prove that the diagram
\[
\xymatrix{{\rm Vect}_G(M)\ar[d]_{{\rm Ch}}\ar[rr]^{{\rm Ch}_\Omega}&&HC^\bullet(\mathcal{A})\\H^\bullet_G(M)\ar[urr]_{\Phi}&&}
\]
commutes when the action is {\em proper}.






\subsection{On the cyclic homology of the crossed product algebra}
Let $\mathsf{A}$ be a commutative associate algebra equipped with a locally convex topology which is acted upon by a Lie group $G$ by means of automorphisms via a homomorphism $\rho:G\to{\rm Aut}(\mathsf{A})$. This action is supposed to be smooth in the following sense, c.f. \cite{Blanc}:
\begin{itemize}
\item[a)] the map $G\times\mathsf{A}\to\mathsf{A}$, given by $(g,a)\mapsto \rho(g)(a)$, is continuous,
\item[b)] for any $a\in\mathsf{A}$, the map $g\mapsto \rho(g)(a)$ is smooth,
\item[c)] $\rho(C)\subset{\rm Aut}(\mathsf{A})$ is equicontinuous for all compact subsets $C\subset G$.
\end{itemize}
For such an algebra $\mathsf{A}$ we can define an action $\rho:\mathfrak{g}\to{\rm Der}(\mathsf{A})$ of the Lie algebra $\mathfrak{g}$ of $G$ by derivations as follows:
\[
\rho(X)(a):=\left.\frac{d}{dt}\right|_{t=0}\rho(\exp(tX))(a),\qquad X\in\mathfrak{g}.
\]
Let $\Omega^\bullet\mathsf{A}=\bigwedge^\bullet_\mathsf{A}\Omega^1\mathsf{A}$ be the DG algebra of K\"ahler differentials of $\mathsf{A}$
 generated by the universal $\mathsf{A}$-module for derivations $\Omega^1\mathsf{A}$. By definition we have ${\rm Der}(\mathsf{A})\cong{\rm 
 Hom}_\mathsf{A}(\Omega^1\mathsf{A},\mathsf{A})$ so we obtain a map $\mathfrak{g}\to {\rm Hom}_\mathsf{A}(\Omega^1\mathsf{A},
 \mathsf{A})$ that we write as $X\mapsto \iota_X,~X\in\mathfrak{g}$ and can extend to $\Omega^\bullet\mathsf{A}$ as a graded derivation of 
 degree $-1$. We can then define a $\mathfrak{g}$-module structure on the whole of $\Omega^\bullet\mathsf{A}$ using Cartan's magic formula 
 $L_X:=d\circ\iota_X+\iota_X\circ d$. Together with the natural $G$-action induced by that on $\mathsf{A}$, this turns
 $\Omega^\bullet\mathsf{A}$ into what is called a $G^\star$-algebra in \cite[\S 2.3]{GS}.

Denote by $G\ltimes\mathsf{A}$ the crossed product algebra associated to the $G$-action on $\mathsf{A}$, and by $CC_\bullet(G\ltimes\mathsf{A})$ its cyclic bicomplex. In this section we shall construct, following Brylinski \cite{brylinski} and Block-Getzler \cite{BG}, chain morphisms
\begin{equation}
\label{cdch}
CC_\bullet(G\ltimes\mathsf{A})\xrightarrow{\Psi_1} {\rm Diag}(CC_{\bullet,\bullet}(G,\mathsf{A}))\xrightarrow{\Psi_2} {\rm Tot}(CC_{\bullet,\bullet}(G,\mathsf{A}))\xrightarrow{\Psi_3} {\rm Tot}(CC_{\bullet,\bullet}(G,\Omega_\mathfrak{g}\mathsf{A})).
\end{equation}
Here $C_{\bullet,\bullet}(G,\mathsf{A})$ forms a double 
complex with horizontal and vertical differentials given by the group homology differential and the (twisted) Hochschild differential of the commutative algebra $\mathsf{A}$, while $C_{\bullet,\bullet}(G,\Omega_\mathfrak{g}\mathsf{A})$ forms a double complex which takes in ideas of Block-Getzlers model for equivariant differential forms (for details, see below).

\subsubsection{The cyclic chain complex}
Because the crossed product algebra $\mathcal{A}:=G\ltimes\mathsf{A}$ is $H$-unital, we can define the space of Hochschild chains
as
\[
C_\bullet(\mathcal{A}):=\begin{cases}\mathcal{A},&\bullet=0,\\\mathcal{A}^+\otimes\mathcal{A}^{\otimes \bullet},&\bullet>0.\end{cases}
\]
Together with the Hochschild differential $b:C_k(\mathcal{A})\to C_{k-1}(\mathcal{A})$ and Connes' cyclic differential 
$B:C_k(\mathcal{A})\to C_{k+1}(\mathcal{A})$ this forms a mixed complex. We denote by $CC_\bullet(\mathsf{A})=\bigoplus_{k=0}^{[n/2]}C_{n-2k}(\mathsf{A})$ the associated complex computing cyclic homology.

\subsubsection{The cylindrical space $L^+(\mathsf{A},G)$} Associated to the action of $G$ on $\mathsf{A}$ is the cylindrical 
vector space $L^+(\mathsf{A},G)$ (c.f.\ \cite{BGJ}) with  $L^+(\mathsf{A},G)_{p,q}:=(C^\infty_c(G)^+)^{\otimes (p+1)}\otimes (\mathsf{A}^+)^{\otimes (q+1)}$, where $C^\infty_c(G)^+$ denotes the smooth convolution algebra on $G$ with a unit adjoined, i.e., the Dirac 
delta function $\delta_e$ at the unit. Because we are using the inductive tensor product we have an isomorphism 
$C^\infty_c(G)\bar{\otimes}C^\infty_c(G)\cong C^\infty_c(G\times G)$ so we can think of elements in $L^+(\mathsf{A},G)_{p,q}$ as 
compactly supported distributional maps from $G^{\times(p+1)}$ to $(\mathsf{A}^+)^{\otimes (q+1)}$ where we only allow 
allow singular distributions of the form $\delta_e(g_i),~i=0,\ldots p$. If we have $F\in L^+(\mathsf{A},G)_{p,q}$ of the form
$F(g_0,\ldots,g_p)=f(g_0,\ldots,g_p)\otimes a_0\otimes\ldots\otimes a_q$ with $f$ a compactly supported distribution on $G^{\times (p+1)}$
of the type as described above and $a_i\in\mathsf{A}^+,~i=0,\ldots,q$, then the maps inducing the cylindrical structure of $L^+(\mathsf{A},G)$ are given by
\begin{align*}
d^h_i(F)(g_0,...,g_p)&=f(g_0,...,g_p)\otimes a_0\otimes\ldots \otimes a_ia_{i+1}\otimes\ldots\otimes a_q&(0\leq i\leq q-1)\\
d^h_q(F)(g_0,...,g_p)&=f(g_0,...,g_p)\otimes (g_0\cdots g_q)^{-1}(a_q)a_0\otimes\ldots\otimes a_{q-1}&\\
s^h_i(F)(g_0,...,g_p)&=f(g_0,...,g_q)\otimes a_0\otimes\cdots\otimes a_i\otimes 1\otimes a_{i+1}\otimes\cdots\otimes a_p&(0\leq i\leq q)\\
t^h(F)(g_0,...,g_p)&=f(g_0,...,g_p)\otimes (g_0\cdots g_q)^{-1}(a_q)\otimes a_0\otimes\cdots\otimes a_{q-1}\\
d^v_i(F)(g_0,...,g_{p-1})&=\int_G F(g_0,...,\gamma,\gamma^{-1}g_i,...,g_{p-1})d\gamma&(0\leq i\leq p-1)\\
d^v_p(F)(g_0,...,g_{p-1})&=\int_G \gamma\cdot F(\gamma^{-1}g_0,g_1,...,g_{p-1},\gamma)d\gamma&\\
s^v_i(F)(g_0,...,g_{p+1})&=\delta(g_{i+1})F(g_0,...,g_i,g_{i+2},...,g_{p+1})&(0\leq i\leq p)\\
t^v(F)(g_0,...,g_p)&=g_0\cdot F(g_1,...,g_p,g_0)&
\end{align*}
From the cylindrical structure, we obtain a double complex
\[C_{p,q}(G,\mathsf{A})\coloneqq C^\infty_c(G)^+\otimes C^\infty_c(G)^{\otimes p}\otimes \mathsf{A}^+\otimes\mathsf{A}^{q+1},\]
which can be thought of as the \textit{normalized} quotient of $L^+(G,\mathsf{A})$ by both the horizontal and vertical degeneracies. On this double complex, the horizontal and vertical differentials are induced by the formulas
\begin{align*}
b^h&=\sum_{i=0}^p (-1)^i d_i^h,&
b^v&=\sum_{i=0}^q (-1)^{i+p}d_i^v.
\end{align*}
Now in the cylindrical space $L^+(G,\mathsf{A})$ we can contract $(b')^h$ and $(b')^v$ using
\begin{align*}
c^h(F)(g_0,...,g_q)&=1\otimes F(g_0,...,g_q)\\
c^v(F)(g_0,...,g_{q+1})&=\delta(g_0)F(g_1,...,g_{q+1})
\end{align*}
and so using general principles given in \cite{gj} we can cook up cyclic differentials $B^h$ and $B^v$ such that $(C_{p,q}(G,\mathsf{A}),b^h+ b^v,B^h+ B^v)$ is a mixed double complex. These cyclic differentials are induced by the formulas
\begin{align*}
B^h&=(1+(-1)^pt^h)c^h\left(\sum_{j=0}^p d_j^h\right) c^h\left(\sum_{i=0}^p (-1)^{ip}(t^h)^i\right),\\
B^v&=(-1)^p(t^h)^{p+1}(1+(-1)^qt^v)c^v\left(\sum_{j=0}^q d_j^v\right)c^v\left(\sum_{i=0}^q(-1)^{iq}(t^v)^i\right).
\end{align*}
From the mixed double complex $C_{\bullet,\bullet}(G,\mathsf{A})$ we make the total cyclic complex
\[\Tot(CC(G,\mathsf{A}))_n=\bigoplus_{\substack{p+q\leq n\\
p+q\equiv n\text{ mod }2}}C_{p,q}(G,\mathsf{A}),\]
where the differential is given by all appropriate applications of the four differentials $b^h$, $b^v$, $B^h$ and $B^v$.
\subsubsection{The morphism $\Psi_1$}
The first morphism in the diagram \eqref{cdch} is the (well-known) map $\Psi_1: C_k((\mathsf{A}\rtimes G)^+)\to C_{k,k}(G,\mathsf{A})$  first  
written down in \cite{brylinski} and given by
\begin{equation*}
\Psi_1(F)(g_0,...,g_n)=((g_0\cdots g_n)^{-1}\otimes\cdots \otimes g_n^{-1})F(g_0,...,g_n)
\end{equation*}
which intertwines the cyclic structure on the Hochschild chains of $\mathsf{A}\rtimes G$ with the cyclic structure on the diagonal 
${\rm Diag}(L^+(\mathsf{A},G)_{\bullet,\bullet})$.

\subsubsection{The morphism $\Psi_2$}
Due to the Eilenberg-Zilber Theorem, the maps $\EZ_{p,q}\colon C_{p+q,p+q}(G,\mathsf{A})\to C_{p,q}(G,\mathsf{A})$ and $\nabla_{p,q}\colon C_{p,q}(G,\mathsf{A})\to C_{p+q,p+q}(G,\mathsf{A})$ given by
\begin{align*}
\EZ_{p,q}&\coloneqq d^h_{p+1}\cdots d^h_{p+q}\cdot (d^v_0)^p,\\
\nabla_{p,q}&\coloneqq \sum_{\sigma\in\Sh(p,q)}(-1)^\sigma s^h_{\sigma(p+q)}\cdots s^h_{\sigma(p+1)}\cdot s^v_{\sigma(p)}\cdots s^v_{\sigma(1)},
\end{align*}
form chain maps $\EZ\colon (\diag(C(G,\mathsf{A})),b)\to (\Tot(C(G,A)),b^h+b^v)$, $\nabla\colon (\Tot(C(G,\mathsf{A})),b^h+b^v)\to(\diag(C(G,\mathsf{A}),b)$ such that
\begin{itemize}
\item $\EZ\circ\nabla=\text{id}$,
\item $\nabla\circ\EZ\simeq\text{id}$ via a homotopy $h: C_{n,n}(G,\mathsf{A})\to C_{n+1,n+1}(G,\mathsf{A})$.
\end{itemize}
By work of Khalkali and Rangipour \cite{kr} we can upgrade this to the cyclic setting by defining a the map $\EZpert: \diag(CC(G,\mathsf{A})\to \Tot(CC(G,\mathsf{A})$ which is defined on $F\in C_{k,k}(G,\mathsf{A})$, seen as an element of $\diag(CC(G,\mathsf{A}))_n$, by
\[\EZpert(F)=\sum_{m=0}^{\frac{n-k}{2}}\EZ((Bh)^m(F)).\]
As is shown in Khalkali and Rangipour, the Homological Pertubation Lemma guarantees us that $\EZpert$ is a quasi-isomorphism.
\subsubsection*{Equivariant differential forms}
Inspired by Getzler's model for equivariant cohomology \cite{getzler}, we employ the equivariant Hochschild-Kostant-Rosenberg-map of Block and Getzler \cite{BG} to obtain a mixed double complex $C_{\bullet,\bullet}(G,\Omega_\mathfrak{g}\mathsf{A})$:
\begin{definition}The mixed double complex $C_{\bullet,\bullet}(G,\Omega_\mathfrak{g}\mathsf{A})$ is defined in degree $(p,q)$ as the subspace
\[C_{p,q}(G,\Omega_\mathfrak{g}\mathsf{A})\subset C^{-\infty}(\mathfrak{g}\times G^{\times q},\Omega^p\mathsf{A})\]
consisting of those distributions that are compactly supported in $G^{\times q}$ and whose singular behaviour is restricted to $\delta(e^X(g_1\cdots g_q)^{-1})$. Associated to these spaces are the four differentials
\begin{align*}
\tilde{b}^h&\colon C_{\bullet,\bullet}(G,\Omega_\mathfrak{g}\mathsf{A})\to C_{\bullet-1,\bullet}(G,\Omega_\mathfrak{g}\mathsf{A}),&\tilde{B}^h&\colon C_{\bullet,\bullet}(G,\Omega_\mathfrak{g}\mathsf{A})\to C_{\bullet+1,\bullet}(G,\Omega_\mathfrak{g}\mathsf{A}),\\
\tilde{b}^v&\colon C_{\bullet,\bullet}(G,\Omega_\mathfrak{g}\mathsf{A})\to C_{\bullet,\bullet-1}(G,\Omega_\mathfrak{g}\mathsf{A}),&\tilde{B}^v&\colon C_{\bullet,\bullet}(G,\Omega_\mathfrak{g}\mathsf{A})\to C_{\bullet,\bullet+1}(G,\Omega_\mathfrak{g}\mathsf{A}),
\end{align*}
given by
\begin{align*}
\tilde{b}^h(F)(X,g_1,...,g_q)&=\iota_X(F(X,g_1,...,g_q)),\\
\tilde{b}^v(F)(X,g_1,...,g_{q-1})&=\int_G(-1)^pF(X,\gamma^{-1},g_1,...,g_{q-1})d\gamma\\
&+\sum_{i=1}^{q-1}(-1)^{i+p}\int_G F(X,g_1,...,\gamma,\gamma^{-1}g_j,...,g_{q-1})d\gamma\\
&+(-1)^{p+q}\int_G \gamma\cdots F(\Ad_{\gamma^{-1}}(X),g_1,...,g_{q-1},\gamma)d\gamma,\\
\tilde{B}^h(F)(X,g_1,...,g_q)&=\int_{\Delta^1}e^{-tX}\cdot d(F(X,g_1,...,g_q))dt\\
\tilde{B}^v(F)(X,g_1,...,g_{q+1})&=\sum_{i=0}^q(-1)^{iq+p}\delta(e^X(g_1\cdots g_{q+1})^{-1})(g_{i+1}\cdots g_{q+1})^{-1}\cdot \\
&\,\,\,\,\,\,\,\,\cdot F(\Ad_{g_{i+1}\cdots g_{q+1}}(X),g_{i+2},...,g_{q+1},g_1,...,g_i).
\end{align*}
\end{definition}
One checks explicitly that with the four operators $\tilde{b}^h$, $\tilde{b}^v$, $\tilde{B}^h$ and $\tilde{B}^v$, $C_{\bullet,\bullet}(G,\Omega_\mathfrak{g}\mathsf{A})$ becomes a double mixed complex.
\subsubsection{The morphism $\Psi_3$}
The last morphism in the diagram \eqref{cdch} is given by an equivariant variant of the Hochschild-Kostant-Rosenberg map as in \cite{BG} with the
difference that here we are not dealing with $G$-invariant chains as in {\em loc. cit.}, but with the equivariant chains from the previous subsection.

We recall the equivariant Hochschild-Kostant-Rosenberg map of Block and Getzler \cite{BG}, which for $X\in\mathfrak{g}$ is given by the map $\HKR_X: A^+\otimes A^{\otimes p}\to \Omega^p\mathsf{A}$ and the formula
\[\HKR_X(a_0\otimes\cdots\otimes a_p)\coloneqq\int_{\Delta^p}a_0d(e^{-t_1X}\cdot a_1)\wedge\cdots\wedge d(e^{-t_pX}\cdot a_p)dt_1\cdots dt_p.\]

With this map we define the morphism $\Psi_3$ on the level of chains as the map $\Psi_3: C_{p,q}(G,\mathsf{A})\to C_{p,q}(G,\Omega_\mathfrak{g}\mathsf{A})$ with the formula
\[\Psi_3(F)(X,g_1,...,g_q)\coloneqq \HKR_X(F(e^X(g_1\cdots g_q)^{-1},g_1,...,g_q)).\]

By an explicit calculation one infers the following:
\begin{lemma}The map $\Psi_3$ satisfies:
\begin{align*}
\tilde{b}^h\circ \Psi_3&=\Psi_3\circ b^h,&\tilde{B}^h\circ\Psi_3&=\Psi_3\circ B^h,\\
\tilde{b}^v\circ\Psi_3&=\Psi_3\circ b^v,&\tilde{B}^v\circ\Psi_3&=\Psi_3\circ B^v.
\end{align*}
\end{lemma}
One can proof this Lemma with an explicit computation.
\begin{corollary}
The map $\Psi_3$ induces a chain map
\[\Psi_3: \Tot(CC(G,\mathsf{A}))\to \Tot(CC(G,\Omega_\mathfrak{g}\mathsf{A})).\]
\end{corollary}

\subsection{On the periodic cyclic homology of the crossed product algebra}
The construction from the previous subsection can also be applied to periodic cyclic homology. The maps $\Psi_1$ and $\Psi_3$ copy directly from the previous section, while the perterbed Eilenberg-Zilber map now takes the form $\EZpert: \diag(CP_{\ast,\ast}(G,\mathsf{A}))\to \Tot(CP_{\ast,\ast}(G,\mathsf{A}))$ given for $F\in C_{k,k}(G,A)\subset \diag(CP_{\ast,\ast}(G,A))_n$ by
\begin{equation*}
\EZpert(F)=\sum_{m\geq 0}\EZ((Bh)^m(F))
\end{equation*}
It should be noted here that at every degree on the codomain this map only has finite contributions, since the grading on $C_{\ast,\ast}$ is bounded from below. However, since an element in the domain gives a priori contributions in infinite degrees we really need the periodic complex to be defined using the direct product.

All in all we obtain a sequence of maps between associacted periodic cyclic complexes:
\begin{equation*}
CP(G\ltimes\mathsf{A})\xrightarrow{\Psi_1} \diag(CP(G,\mathsf{A}))\xrightarrow{\EZpert} \Tot(CP(G,\mathsf{A}))\xrightarrow{\Psi_3} \Tot(CP(G,\Omega_\mathfrak{g}\mathsf{A}))
\end{equation*}
and the total composition will be denoted by $\Psi:CP(\mathsf{A}\rtimes G)\to \Tot(CP(G,\Omega_\mathfrak{g}\mathsf{A}))$.
\subsection{Pairing with equivariant cohomology}  For this section, let $\AA$ be a smooth unital $G$-algebra containing $\mathsf{A}$ as an ideal, such that $\Omega\mathsf{A}$ is an ideal of $\Omega\AA$. The main example to have in mind is $(\mathsf{A},\AA)=(C^\infty_c(M),C^\infty(M))$ when $M$ is a manifold with a right $G$-action.

We want to pair our presentation of cyclic homology of the crossed product algebra $G\ltimes\mathsf{A}$ with the equivariant cohomology of $\AA$. We recall the model for equivariant cohomology obtained in \cite{getzler} which actually was the inspiration for our presentation of cyclic homology. 

The cochain complex of Getzler is given by
\begin{equation*}
C^{p,q}(G,\Omega_\mathfrak{g}\AA)=\left\{F\in C^\infty(\mathfrak{g}\times G^{\times q},\Omega^p\AA):\begin{matrix}
F\text{ is polynomial in }\mathfrak{g}\\
F(X,g_1,...,g_q)=0\text{ if either one of }g_1,...,g_q\text{ equals 1}
\end{matrix}\right\}
\end{equation*}
endowed with 4 differentials
\begin{align*}
\iota\colon\,&C^{p,q}(G,\Omega_\mathfrak{g}\AA)\to C^{p-1,q}(G,\Omega_\mathfrak{g}\AA)\\
\overline{\iota}\colon\,&C^{p,q}(G,\Omega_\mathfrak{g}\AA)\to C^{p,q-1}(G,\Omega_\mathfrak{g}\AA)\\
d\colon\,&C^{p,q}(G,\Omega_\mathfrak{g}\AA)\to C^{p+1,q}(G,\Omega_\mathfrak{g}\AA)\\
\overline{d}\colon\,&C^{p,q}(G,\Omega_\mathfrak{g}\AA)\to C^{p,q+1}(G,\Omega_\mathfrak{g}\AA)
\end{align*}
given by
\begin{align*}
(\iota F)(X,g_1,...,g_q)=&\,\,(-1)^q\iota_X(F(X,g_1,...,g_k))\\
(\overline{\iota} F)(X,g_1,...,g_{q-1})=&\,\,\sum_{i=0}^{q-1}(-1)^i\left.\frac{d}{dt}\right|_{t=0}F(X,g_1,...,g_i,e^{t\Ad_{g_{i+1}\cdots g_{q-1}}(X)},g_{i+1},...,g_{q-1})\\
(dF)(X,g_1,...,g_q)=&\,\,(-1)^qd(F(X,g_1,...,g_q))\\
(\overline{d}F)(X,g_1,...,g_{q+1})=&\,\,F(X,g_2,...,g_{q+1})\\
&+\sum_{i=1}^q(-1)^iF(X,g_1,...,g_ig_{i+1},...,g_{q+1})\\
&+(-1)^{q+1}g^{-1}_{q+1}\cdot F(\Ad_{g_{q+1}}(X),g_1,...,g_q)
\end{align*}
For the case $\AA=C^\infty(M)$ this complex calculates the equivariant cohomology $H_G(M)$ of $M$ \cite[1.2.3]{getzler}.

In Getzler's work the grading of this complex is the sum of its degree as a group cochain ($q$), its degree as an element of $\Omega\AA$ ($p$) and twice the polynomial degree (not denoted above). With this grading $\iota+\overline{\iota}+d+\overline{d}$ is a cohomological differential. For our deliberations it will be more natural to disregard the polynomial degree and see $C^{\bullet,\bullet}(G,\Omega_\mathfrak{g}\AA)$ as a mixed double cochain complex with differentials $(d+\overline{d},\iota+\overline{\iota})$.

To pair this complex with our complex $C_{\ast,\ast}(G,\Omega_\mathfrak{g}\mathsf{A})$, let $\fint: \Omega^n\mathsf{A}\to\mathbb{R}$ be a closed, graded and $G$-invariant trace (i.e. $\fint d\alpha=0$ and $\fint g\cdot \alpha=\fint\alpha$). We will work with the implicit assumption that $\Omega^{>n}\mathsf{A}=\{0\}$. Using this trace we can write down a pairing $\langle-,-\rangle$ between $C^{n-p,q}(G,\Omega_\mathfrak{g}\AA)$ and $C_{p,q}(G,\Omega_\mathfrak{g}\mathsf{A})$, given for $\alpha_{n-p,q}\in C^{n-p,q}(G,\Omega_\mathfrak{g}\AA)$ and $\beta_{p,q}\in C_{p,q}(G,\Omega_\mathfrak{g}\mathsf{A})$ by
\begin{equation*}
\langle \alpha_{n-p,q},\beta_{p,q}\rangle=(-1)^{p(n+q)+\frac{1}{2}p(p+1)}\int_{G^{\times q}}\fint \alpha_{n-p,q}(0,g_1,...,g_q)\wedge\beta_{p,q}(0,g_1,...,g_q)dg_1\cdots dg_q
\end{equation*}
\begin{lemma}
The following identities hold for all $\alpha$ and $\beta$
\begin{align*}
\langle \iota \alpha_{n-p+1,q},\beta_{p,q}\rangle&=\langle \alpha_{n-p+1,q},\tilde{b^h}\beta_{p,q}\rangle=0\\
\langle d\alpha_{n-p-1,q},\beta_{p,q}\rangle&=\langle \alpha_{n-p-1,q},\tilde{B^h}\beta_{p,q}\rangle\\
\langle \overline{d}\alpha_{n-p,q-1},\beta_{p,q}\rangle&=\langle \alpha_{n-p,q-1},\tilde{b^v}\beta_{p,q}\rangle
\end{align*}
If furthermore $\alpha_{n-p,q+1}$ satisfies that $\alpha_{n-p,q+1}(0,g_1,...,g_{q+1})=0$ when $g_1\cdots g_{q+1}=1$, the following also holds for all $\beta_{p,q}$:
\begin{equation*}
\langle\overline{\iota}\alpha_{n-p,q+1},\beta_{p,q}\rangle=\langle\alpha_{n-p,q+1},\tilde{B^v}\beta_{p,q}\rangle=0
\end{equation*}
\end{lemma}
\begin{proof}
For the first equation we note that $(\iota\alpha_{n-p+1,q})(0,g_1,...,g_q)=0$ since $\iota_0\omega=0$ for any $\omega\in \Omega\mathsf{A}$. Similarly $(\tilde{b^h}\beta_{p,q})(0,g_1,...,g_q)=0$.

For the second line note that for $\omega_1\in\Omega^{n-p-1}\mathsf{A}$ and $\omega_2\in\Omega^p\mathsf{A}$ we have
\begin{equation*}
d(\omega_1\wedge \omega_2)=d\omega_1\wedge \omega_2+(-1)^{n-p-1}\omega_1\wedge d\omega_2
\end{equation*}
Since $\fint$ vanishes on exact forms, we obtain
\begin{equation*}
\fint d\omega_1\wedge \omega_2=(-1)^{n-p}\fint\omega_1\wedge d\omega_2
\end{equation*}
Furthermore we have
\begin{align*}
(\tilde{B^h}\beta_{p,q})(0,g_1,...,g_q)&=\int_{\Delta^1}e^{-t\cdot 0}\cdot d(\beta_{p,q}(0,g_1,...,g_q))dt\\
&=\int_{\Delta^1}d(\beta_{p,q}(0,g_1,...,g_q))dt\\
&=d(\beta_{p,q}(0,g_1,...,g_q))
\end{align*}
This results in
\begin{align*}
\langle d\alpha_{n-p-1,q},\beta_{p,q}\rangle&=\begin{multlined}[t](-1)^{p(n+q)+\frac{1}{2}p(p+1)+q}\int_{G^q}\fint d(\alpha_{n-p-1,q}(0,g_1,...,g_q))\wedge \\
\wedge\beta_{p,q}(0,g_1,...,g_q)dg_1\cdots dg_q\end{multlined}\\
&=\begin{multlined}[t](-1)^{(p-1)(n+q)+\frac{1}{2}p(p-1)}\int_{G^q}\fint \alpha_{n-p-1,q}(0,g_1,...,g_q)\wedge\\
\wedge d(\beta_{p,q}(0,g_1,...,g_q))dg_1\cdots dg_q\end{multlined}\\
&=\begin{multlined}[t](-1)^{(p-1)(n+q)+\frac{1}{2}p(p-1)}\int_{G^q}\fint \alpha_{n-p-1,q}(0,g_1,...,g_q)\wedge\\
\wedge(\tilde{B^h}\beta_{p,q})(0,g_1,...,g_q)dg_1\cdots dg_q
\end{multlined}\\
&=\langle \alpha_{n-p-1,q},\tilde{B^h}\beta_{p,q}\rangle.
\end{align*}
Then for the third part of the Lemma we do an explicit calculation:
\begin{align*}
\langle \overline{d}\alpha_{n-p,q-1},\beta_{p,q}\rangle&=\begin{multlined}[t](-1)^{p(n+q)+\frac{1}{2}p(p+1)}\int_{G^{\times q}}\fint \alpha_{n-p,q-1}(0,g_2,...,g_q)\wedge\\
\wedge \beta_{p,q}(0,g_1,...,g_q)dg_1\cdots dg_q\end{multlined}\\
&+\begin{multlined}[t](-1)^{p(n+q)+\frac{1}{2}p(p+1)}\sum_{i=1}^{q-1}(-1)^i\int_{G^{\times q}}\fint\alpha_{n-p,q-1}(0,g_1,...,g_ig_{i+1},...,g_q)\wedge\\
\wedge \beta_{p,q}(0,g_1,...,g_q)dg_1\cdots dg_q\end{multlined}\\
&+\begin{multlined}[t](-1)^{p(n+q)+\frac{1}{2}p(p+1)}(-1)^q\int_{G^{\times q}}\fint (g_q^{-1}\cdot \alpha_{n-p,q-1}(0,g_1,...,g_{q-1}))\wedge \\
\wedge\beta_{p,q}(0,g_1,...,g_q)dg_1\cdots dg_q.\end{multlined}
\end{align*}
Now using a few change of variables for the integrals of $G^{\times q}$ (and using the fact that $G$ is unimodular for the first line) and using that $\fint$ is $G$-invariant, we obtain:\begin{align*}
\langle d\overline{\alpha}_{n-p,q-1},\beta_{p,q}\rangle=&\begin{multlined}[t](-1)^{p(n+q-1)+\frac{1}{2}p(p+1)}(-1)^p\int_{G^{\times (q-1)}}\int_G\fint\alpha_{n-p,q-1}(0,g_1,...,g_{q-1})\wedge\\
\wedge \beta_{p,q}(0,\gamma^{-1},g_1,...,g_{q-1})d\gamma dg_1\cdots dg_q\end{multlined}\\
&+\begin{multlined}[t](-1)^{p(n+q-1)+\frac{1}{2}p(p+1)}\sum_{i=1}^{q-1}(-1)^{p+i}\int_{G^{\times(q-1)}}\int_G\fint\alpha_{n-p,q-1}(0,g_1,...,g_{q-1})\wedge\\
\wedge\beta_{p,q}(0,g_1,...,\gamma,\gamma^{-1}g_i,...,g_{q-1})dg_1\cdots dg_{q-1}\end{multlined}\\
&+\begin{multlined}[t](-1)^{p(n+q-1)+\frac{1}{2}p(p+1)}(-1)^{p+q}\int_{G^{\times (q-1)}}\int_G\fint\alpha_{n-p,q-1}(0,g_1,...,g_{q-1})\wedge\\
\wedge (\gamma\cdot \beta_{p,q}(0,g_1,...,g_{q-1},\gamma))d\gamma dg_1\cdots dg_{q-1}\end{multlined}\\
=&\begin{multlined}[t](-1)^{p(n+q-1)+\frac{1}{2}p(p+1)}\int_{G^{\times (q-1)}}\fint \alpha_{n-p,q-1}(0,g_1,...,g_{q-1})\wedge\\
\wedge (\tilde{b^v}\beta_{p,q})(0,g_1,...,g_{q-1})dg_1\cdots dg_{q-1}
\end{multlined}\\
=&\langle \alpha_{n-p,q-1},\tilde{b^v}\beta_{p,q}\rangle.
\end{align*}
Then for the last line, note that $(\overline{\iota}\alpha_{n-p,q+1})(0,g_1,...,g_q)=0$, while if one were to write down $\tilde{B^v}$ and calculate the pairing, using that $\alpha_{n-p,q+1}(0,1,g_2,...,g_{q+1})=0$ one obtains
\begin{multline*}
\langle \alpha_{n-p,q+1},\tilde{B^v}\beta_{p,q}\rangle=\int_{G^{\times q}}\fint\sum_{i=0}^q(-1)^{iq}(g_{q-i+1}\cdots g_q)\cdot\\
\cdot\alpha_{n-p,q+1}(0,g_{q-i+1},...,g_q,(g_1\cdots g_q)^{-1},g_1,...,g_{q-i})\wedge\\
\wedge \beta_{p,q}(0,g_1,...,g_q)dg_1\cdots dg_q
\end{multline*}
which vanishes since $g_{q-i+1}\cdots g_q\cdot (g_1\cdots g_q)^{-1}\cdot g_1\cdots g_{q-i}=1$ and so we can use the extra assumption on $\alpha_{n-p,q+1}$.
\end{proof}
From the last part of the previous lemma we arrive at the following:
\begin{definition}
An element $\alpha_{p,q}\in C^{p,q}(G,\Omega_\mathfrak{g}\AA)$ is called \textit{cyclically normalized} if \[\alpha_{p,q}(0,g_1,...,g_q)=0\] when $g_1\cdots g_q=1$.
\end{definition}
In the end we would like to have an actual chain map from $C^{\bullet,\bullet}(G,\Omega_\mathfrak{g}\AA)$ to the periodic cyclic complex of $G\ltimes\mathsf{A}$. However, as is evident from the degrees that get paired, this is not possible if we impose the standard grading on $C^{\bullet,\bullet}(G,\Omega_\mathfrak{g}\AA)$.
To counteract this, we repack Getzler's complex akin to the periodic cyclic complex, thereby also justifying the fact that we disregard the polynomial degree in Getzler's complex.

\begin{definition}The complex $CP^\bullet(G,\Omega_\mathfrak{g}\AA)$ is given by
\begin{equation*}
CP^k(G,\Omega_\mathfrak{g}\mathsf{A})=\bigoplus_{p+q\equiv n+k\text{ mod }2}C^{p,q}(G,\Omega_\mathfrak{g}\mathsf{A})
\end{equation*}
with differential given by $\iota+\overline{\iota}+d+\overline{d}$.
\end{definition}
As $CP^\bullet(G,\Omega_\mathfrak{g}\AA)$ is just the direct sum of all shifts of ${\rm Tot}(C^{\bullet,\bullet}(G,\Omega_\mathfrak{g}\AA))$ by an even number of degrees, the following is clear:
\begin{lemma}
The complex $CP^\bullet(G,\Omega_\mathfrak{g}\AA)$ calculates the following cohomology
\begin{equation*}
H^k(CP^\bullet(G,\Omega_\mathfrak{g}\AA))=\bigoplus_{p\equiv n+k\text{ mod }2}H^p({\rm Tot}(C^{\bullet,\bullet}(G,\Omega_\mathfrak{g}\AA)))
\end{equation*}
\end{lemma}
The previous discussion can then be summarized by the following proposition.
\begin{proposition}
The map $\Phi: CP^\bullet(G,\Omega_\mathfrak{g}\AA)\to (CP(G,\Omega_\mathfrak{g}\mathsf{A})^\times)^\bullet$ given by
\begin{equation*}
\Phi(\sum_{p',q'}\alpha_{p',q'})(\beta_{p,q})=\langle\alpha_{n-p,q},\beta_{p,q}\rangle
\end{equation*}
is compatible with differentials when restricted to cyclically normalized cochains $\sum_{p',q'}\alpha_{p',q'}$.
\end{proposition}
Combining with the map $\Psi$ we arrive at the first main result of this section, which is a corollary of the previous proposition and the fact that $\Psi$ is a chain map.
\begin{corollary}\label{almostchainmap}
The map $c: CP^\bullet(G,\Omega_\mathfrak{g}\AA)\to CP^\bullet(G\ltimes \mathsf{A})$ given by
\begin{equation*}
c(\sum_{p,q}\alpha_{p,q})(a_0,...,a_k)=\Phi(\sum_{p,q}\alpha_{p,q})(\Psi(a_0\otimes\cdots\otimes a_k))
\end{equation*}
is compatible with differentials when restricted to cyclically normalized cochains $\sum_{p,q}\alpha_{p,q}$.
\end{corollary}
\subsection{The case $\AA=C^\infty(M)$}
As we have seen in the previous part, the fact that not every element of $C^{\bullet,\bullet}(G,\Omega_\mathfrak{g}\AA)$ is cyclically normalized means that we cannot write down a chain map on the nose. When we return to the case of a manifold $M$ with a right $G$-action we can explicitely circumvent this problem to obtain an actual map from equivariant cohomology to periodic cyclic cohomology, in a way we describe now.

So, in this section $(\mathsf{A},\AA)=(C^\infty_c(M),C^\infty(M))$, where $M$ is a smooth oriented manifold of dimension $n$ with a right orientation-preserving $G$-action. Note that in this case integration of top-forms induces a closed, graded, $G$-invariant trace $\int_M$ on $\Omega^n\mathsf{A}$.

The clue is to quasi-isomorphically embed Getzler's complex into a complex with a cyclic structure. To wit, we recall that in \cite{getzler} the complex $C^{\bullet,\bullet}(G,\Omega_\mathfrak{g}\AA)$ is obtained by constructing the reduced cobar resolution $\Omega^\bullet_{\rm red}(\mathbb{C},\Omega^\bullet(G),\Omega\AA)$ and then applying a quasi-isomorphism $\mathcal{J}:\Omega^\bullet_{\rm red}(\mathbb{C},\Omega^\bullet(G),\Omega\AA)\to C^{\bullet,\bullet}(G,\Omega_\mathfrak{g}\AA)$.

The reduced cobar resolution is defined by $\Omega^\bullet_{\rm red}(\mathbb{C},\Omega^\bullet(G),\Omega\AA)=\overline{\Omega^\bullet(G)}\otimes\Omega\AA$ where $\overline{\Omega^\bullet(G)}$ denotes the kernel of the counit $\Omega^\bullet(G)\to\mathbb{C}$ which is evaluation at the identity. The quasi-isomorphism $\mathcal{J}$ is then given by the formula
\[
\mathcal{J}(\omega_1\otimes\ldots\otimes\omega_k\otimes\gamma)(g_1,\ldots,g_\ell,X)=(-1)^l\sum_{\sigma\in {\rm Sh}_{\ell,k-\ell}}\left(\prod_{i=1}^\ell\omega_{\sigma(i)}(g_i)\right)\left(\prod_{j=\ell+1}^k\omega_{\sigma(j)}(e)(X_j)\right)
\gamma.
\]
where $X_j=\Ad_{g_i\cdots g_l}X$ with $i\leq l$ minimal such that $\sigma(j)<\sigma(i)$ (and $X_j=X$ if such an $i$ does not exist). In this formula the terms $\omega(g)$ only contribute when $\omega$ is a zero-form on $G$, and the terms $\omega(e)(X)$ only contribute when $\omega$ is a one-form on $G$.

We remark that using the kernel of the counit $\Omega^\bullet(G)\to\mathbb{C}$ in the definition of the reduced cobar-resolution corresponds directly to the fact that in Getzler's complex $C^{\bullet,\bullet}(G,\Omega_\mathfrak{g}\AA)$ the chains satisfy that $\alpha(X,g_1,...,g_q)=0$ whenever $e\in\{g_1,...,g_q\}$.

However, we can equivalently use the unreduced cobar resolution $\Omega^\bullet(\mathbb{C}
,\Omega^\bullet(G),\Omega\AA)=\Omega^\bullet(G)\otimes\Omega\AA$, to which the formula of $\mathcal{J}$ naturally extends. The result will be a complex $\tilde{C}^{\bullet,\bullet}(G,\Omega_\mathfrak{g}\AA)$ which is a quasi-isomorphic extension of Getzlers complex $C^{\bullet,\bullet}(G,\Omega_\mathfrak{g}\AA)$ given by
\begin{equation*}
\tilde{C}^{p,q}(G,\Omega_\mathfrak{g}\AA)=\left\{F\in C^\infty(\mathfrak{g}\times G^{\times q},\Omega^p\AA):\begin{matrix}
F\text{ is polynomial in }\mathfrak{g}
\end{matrix}\right\}
\end{equation*}

Now for the case $\AA=C^\infty(M)$ we note that the unreduced cobar resolution is simply $\Omega^\bullet(G^{\times \ast}\times M)$ with the differential the sum of the deRham-differential and the differential coming from the underlying simplicial structure on $G^{\times\bullet}\times M$ given as usual by face operator $\partial_i:G^{\times k}\times M\to G^{\times(k-1)}\times M$ defined as
\[
\partial_i(g_1,\ldots,g_k,x)=\begin{cases} (g_1,\ldots,g_ig_{i+1},\ldots,g_k,x),& 0\leq i\leq k-1,\\
(g_1,\ldots,g_{k-1},g_kx),& i=k.\end{cases}
\]
This simplicial space also has a cyclic structure given by
\[
t(g_1,\ldots,g_k,x)= ((g_1\cdots g_k)^{-1},g_1,\ldots,g_{k-1},g_kx)
\]
The induced map by pull-back on differential forms gives the unreduced cobar resolution
the structure of a DG cocyclic vector space.

When computing the cyclic cohomology, we can, as usual restrict ourselves to the cyclic subcomplex, in this case given 
by
\[
\mathsf{\Omega}^\bullet(\C,\Omega G,\Omega M)_\lambda:=\bigoplus_q\{\alpha\in \Omega(G^{\times q}\times M),~t^*\alpha=(-1)^q\alpha\}.
\]
We call forms satisfying $t^*\alpha=(-1)^q\alpha$ {\em cyclic differential forms}. It is then clear that the cohomology of $\tilde{C}^{\bullet,\bullet}(G,\Omega_\mathfrak{g}M)$ is exhausted by classes obtained by applying $\mathcal{J}$ to such cyclic differential forms.

Next we note that the pairing between $C^{\bullet,\bullet}(G,\Omega_\mathfrak{g}\AA)$ and $C_{\ast,\ast}(G,\Omega_\mathfrak{g}\mathsf{A})$ can naturally be extended to $\tilde{C}^{\ast,\ast}(G,\Omega_\mathfrak{g}\AA)$, and we will show that the cyclic differential forms actually induce cyclically normalized elements. First we have the following observation that stems from an elementary calculation:
\begin{lemma}
The image of a cyclic differential form $\alpha\in \Omega(G^{\times \bullet}\times M)_\lambda$ satisfies
\[
\mathcal{J}(\alpha)(0,g_1,\ldots,g_q)=(-1)^q(g_q^{-1})^\ast(\mathcal{J}(\alpha)(0,(g_1\cdots g_q)^{-1},g_1,\ldots,g_{q-1}))
\]
\end{lemma}
Then we can again do a computation with the differential $\tilde{B^v}$ to obtain the following:
\begin{lemma}
Let $\alpha_{p,q}\in\tilde{C}^{p,q}(G,\Omega_\mathfrak{g}\AA)$ be such that
\begin{equation*}
\alpha_{p,q}(0,g_1,...,g_q)=(-1)^q(g_q^{-1})^\ast\alpha(0,(g_1\cdots g_q)^{-1},g_1,...,g_{q-1})
\end{equation*}
then for all $\beta_{n-p,q-1}\in C_{n-p,q-1}(G,\Omega_\mathfrak{g}\mathsf{A})$ we have
\begin{equation*}
\langle \alpha_{p,q},\tilde{B^v}\beta_{n-p,q-1}\rangle=0
\end{equation*}
\end{lemma}
This is clear from the definition of $\tilde{B^v}$ and the defining equation for a cyclic differential form.
\begin{corollary}
The cohomological pairing between $\Omega^\bullet(G^{\times \bullet}\times M)_\lambda$ and $C_{\bullet,\bullet}(G,\Omega_\mathfrak{g}(M))$, combined with the construction the map $\Psi: CP_\bullet(G\rtimes C^\infty_c(M))\to CP_\bullet(\Tot(C(G,\Omega_\mathfrak{g}M))$, induces a map $c: H_G^{\rm{ev/odd}}(M)\to HP^{\rm{ev/odd}+n}(G\ltimes C^\infty(M))$ which, applied to such classes in $H^{\rm{ev/odd}}_G(M)$ induced by cyclically normalized elements of $C^{\bullet,\bullet}(G,\Omega_\mathfrak{g}(M))$, agrees with the map $c$ from \ref{almostchainmap}.
\end{corollary}
\subsection{Compatibility with the Chern character}
In \cite{getzler} the equivariant Chern character of an equivariant vector bundle $E\to M$ with connection $\Nabla$ (not necessarily invariant) was determined by an explicit cyclic cocycle $\Ch_G(E,\Nabla)$ of degree $n$ in $\Tot(C^{\ast,\ast}(G,\Omega_\mathfrak{g}(M)))$. For our pairing we are only interested in $\Ch_G(E,\Nabla)^0_q\in C^\infty(G^{\times q},\Omega(M))$ given by \[\Ch_G(E,\Nabla)^0_q(g_1,...,g_q)\coloneqq\Ch_G(E,\Nabla)(0,g_1,...,g_q)\] and we describe these functions.

They are given by the same formulae as the Chern-Simons forms of Bott when $G$ is discrete. To wit, let $\Nabla_0,...,\Nabla_q$ be connections of $E$. The Chern-Simons form is then defined by

\begin{align*}
\text{cs}(\Nabla_0,...,\Nabla_q)&=\int_{\Delta^q}\Tr\left(\exp\left(-dt_1\wedge(\Nabla_q-\Nabla_0)+\sum_{i=1}^{q-1} (dt_i-dt_{i+1})\wedge (\Nabla_q-\Nabla_i)+F(t)\right)\right)\\
&=\int_{\Delta^q}\Tr\left(\exp\left(\sum_{i=1}^qdt_i\wedge(\Nabla_{i-1}-\Nabla_i)+F(t)\right)\right)
\end{align*}
where $F(t)$ is the curvature of the connection
\begin{equation*}
\Nabla(t)=\Nabla_q-t_1(\Nabla_q-\Nabla_0)+\sum_{i=1}^{q-1}(t_i-t_{i+1})(\Nabla_q-\Nabla_i)
\end{equation*}

Then closely dissecting the construction of Getzler, we recognize that
\begin{equation*}
\text{Ch}_G(E,\Nabla)^0_{p,q}(g_1,...,g_q)=(-1)^{p+q}\text{cs}(\gamma_1^\ast\Nabla,...,\gamma_q^\ast\Nabla,\Nabla)
\end{equation*}
where $\gamma_i=g_i\cdots g_q$.

Investigating the specific from of the Chern-Simons forms then gives the following Lemma.
\begin{lemma}
The functions $\Ch_G(E,\Nabla)^0_q$ satisfy the following:
\begin{itemize}
\item[i)] $\Ch_G(E,\Nabla)^0_q(g_1,...,g_q)\in\bigoplus_{i=q}^n\Omega^i\mathsf{A}$
\item[ii)] $\Ch_G(E,\Nabla)^0_q(g_1,...,g_q)=0$ if either one of the $g_i$'s or their product $g_1\cdots g_q$ is the unit of the group.
\end{itemize}
\end{lemma}
\begin{proof}
For i) note that the only terms under the exponent in the Chern-Simons form that persist in $\Omega\mathsf{A}$ after taking the integral over $\Delta^q$ are the terms where all the parts $dt_i\wedge (\Nabla_{i-1}-\Nabla_i)$ appear exactly one. This means in particular that the term that results from $\int_{\Delta^q}$ is a form on $M$ of degree at least $q$.

For ii) note that $g_i=1$ for $i=1,...,q-1$ correspond to $\gamma_i^\ast\Nabla=\gamma_{i+1}^\ast\Nabla$, $g_q=1$ corresponds to $\gamma_q^\ast\Nabla=\Nabla$ and $g_1\cdots g_q=1$ corresponds to $\gamma_1^\ast\Nabla=\Nabla$. To obtain the statement of the Lemma, we now argue that $\cs(\Nabla_0,...,\Nabla_q)=0$ if either $\Nabla_{i-1}=\Nabla_i$ for some $i=1,...,q$ or $\Nabla_0=\Nabla_q$. Indeed, if one of these holds we look at the form on $\Delta^q\times M$ of which we take the exponent and look at the forms on $\Delta^q$ which are part of it. If we have one the equalities of the $\Nabla_i$'s we see that only $q-1$ different one forms on $\Delta^q$ remain (for the case $\Nabla_{i-1}=\Nabla_i$ we look at the second way of writing $\cs(\Nabla_0,...,\Nabla_q)$ and only the terms $\{dt_j\wedge(\Nabla_{j-1}-\Nabla_j)\}_{j\neq i}$ remain, while for the case $\Nabla_q=\Nabla_0$ we use the first way of writing $\cs(\Nabla_0,...,\Nabla_q)$ to see that only $\{(dt_i-dt_{i+1})\wedge (\Nabla_q-\Nabla_i)\}_{i=1,...,q-1}$ remain), then the only way to obtain a $q$-form on $\Delta^q$ is to take a $q$-fold wedge product with at least one repeating one-form, and hence all the terms vanish.
\end{proof}

\begin{remark}
The second part of previous Lemma show that the periodic cyclic class in $HP^n(G\ltimes C^\infty(M))$ induced by the equivariant Chern character $[\Ch_G(E,\nabla)]\in H^{\rm{even}}_G(M)$ is represented by the explicit periodic cyclic cochain $c(\Ch_G(E,\nabla))$ obtained by \ref{almostchainmap}. Indeed, from it one concludes that $\Ch_G(E,\nabla)$ is a cyclically normalized cochain.
\end{remark}
\begin{remark}
The first part of the previous Lemma implies that $c(\Ch_G(E,\nabla))$ is actually a cyclic cochain of degree $n$, not just a periodic cyclic cochain.
\end{remark}

We now come to the main result of this paper, namely that the two equivariant Chern characters we constructed agree in periodic cyclic cohomology when the action is proper:

\begin{theorem}\label{chernsagree}
If the action of $G$ on $M$ is proper, the Chern character ${\rm Ch}_{\Omega_{E,\nabla}}$ of \cref{ncchern} and the periodic cyclic cocycle $c(\Ch_G(E,\nabla))$ obtained from plugging in Getzler's Chern character into Theorem \ref{almostchainmap} induce the same periodic cyclic cohomology classes in $HP^n(G\ltimes C^\infty(M))$.
\end{theorem}

The gist of the proof is that when the action is proper, there is a $G$-invariant connection on $E$.

\begin{remark}
If is easy to see that the periodic cyclic classes induced by $\Ch_{\Omega_{E,\nabla}}$ and $c(\Ch_G(E,\nabla))$ respectively are both independent of the choice of the connection. For $\Ch_{\Omega_{E,\nabla}}$ we already made this remark in \cref{ncchern}, while for $c(\Ch_G(E,\nabla))$ this follows from the fact that our map $c$ induces a map on cohomology and $[\Ch_G(E,\nabla)]\in H^\ev_G(M)$ is independent of the connection $\nabla$.
\end{remark}

\begin{proof}[Proof of \ref{chernsagree}]
Since the action on $M$ is proper we can find a $G$-invariant connection $\nabla$ on $E$. This has two implications:
\begin{itemize}
\item In $\Omega_{E,\nabla}$ we have $\delta\equiv 0$, in particular $F(\nabla)$ is also $G$-invariant
\item Of the Chern character $\Ch_G(E,\nabla)$ in Getzler's complex, the for us relevant part $\Ch_G(E,\nabla)^0$ only lives in degrees $(2k,0)$ where it is given by $\frac{1}{k!}{\rm tr}(F(\nabla)^{\wedge k})$
\end{itemize}
We can now write down explicitly the two cycles $\Ch_{\Omega_{E,\nabla}}$ and $c(\Ch_G(E,\nabla))$. Starting with the former, we can simplify the formula of \eqref{classoftwistedCuDGA} from before since $G$ is unimodular and obtain that
\begin{multline*}
\Ch_{\Omega_{E,\nabla}}^k(a_0,...,a_k)=\frac{(-1)^{\frac{n-k}{2}}}{k!\left(\frac{n-k}{2}\right)!}\int_M{\rm tr}_E(F(\nabla)^{\wedge\frac{n-k}{2}})\\
\int_{G^k}a_0((h_1\cdots h_k)^{-1})((h_1\cdots h_k)^{-1})^\ast da_1(h_1)\wedge\cdots\\
\cdots (h_k^{-1})^\ast da_k(h_k)dh_1\cdots dh_k
\end{multline*}
This follows from the formula earlier, since $a_0$ and the $a_i$'s are scalar forms inside $\Omega(M,{\rm End}(E))$ and hence commute with the $2$-form $F(\nabla)$, the fact that in the end the summand does not depend on $(i_0,...,i_k)$ and the fact that we sum over an index set of size $\begin{pmatrix}
\frac{n+k}{2}\\k
\end{pmatrix}$.

On the other hand, we can carefully write down the formula for $c(\Ch_G(E,\nabla))$. Starting with $\alpha\in CP^0(G,\Omega_\mathfrak{g}(M))$, we see that the formula for $c(\alpha)$ becomes
\begin{equation}\label{classofpairingEZetc}
c(\alpha)(a_0,...,a_k)=\sum_{m=0}^{\frac{n-k}{2}}\sum_{p=0}^{k+2m}\langle \alpha_{n-p,k+2m-p},\Psi_2\EZ_{p,k+2m-p}((Bh)^m\Psi_1(a_0\otimes\cdots\otimes a_k))\rangle
\end{equation}

For $\alpha=\Ch_G(E,\nabla)$ we see that only the terms of $m=0$ contribute. Indeed, starting with $a_0,...,a_k\in C^\infty_c(G,C^\infty_c(M))$, we obtain that
\begin{multline*}
(\Psi_2\EZ_{k,0}\Psi_1(a_0\otimes\cdots\otimes a_k))(0)=\frac{1}{k!}\int_{G^{\times k}}(h_1\cdots h_k)^\ast a_0((h_1\cdots h_k)^{-1})da_1(h_1)\wedge h_1^\ast da_2(h_2)\wedge\cdots\\
\cdots\wedge (h_1\cdots h_{k-1})^\ast da_k(h_k)dh_1\cdots dh_k
\end{multline*}
here the factor $\frac{1}{k!}$ comes from the fact that plugging in $X=0$ yields and integrand over $\Delta^k$ independent of $t_1,...,t_k$, so we get a factor ${\rm vol}(\Delta^k)$.

Then pairing with $\Ch_G(E,\nabla)$, the $m=0$ terms of Equation \ref{classofpairingEZetc} equal
\begin{multline*}
\frac{(-1)^{\frac{n-k}{2}}}{k!\left(\frac{n-k}{2}\right)!}\int_M{\rm tr}_E(F(\nabla)^{\frac{n-k}{2}})\int_{G^{\times k}}(h_1\cdots h_k)^\ast a_0((h_1\cdots h_k)^{-1})da_1(h_1)\wedge h_1^\ast da_2(h_2)\wedge\cdots\\
\cdots\wedge (h_1\cdots h_{k-1})^\ast a_k(h_k)dh_1\cdots dh_k
\end{multline*}
Then using the fact that $\int_M$ is invariant under the $G$-action in general, and since $F(\nabla)$ is $G$-invariant in this specific case, this can we rewritten as
\begin{multline*}
\frac{(-1)^{\frac{n-k}{2}}}{k!\left(\frac{n-k}{2}\right)!}\int_M{\rm tr}_E(F(\nabla)^{\frac{n-k}{2}})\int_{G^{\times k}}a_0((h_1\cdots h_k)^{-1})((h_1\cdots h_k)^{-1})^\ast da_1(h_1)\wedge\cdots\\
\cdots\wedge (h_k^{-1})^\ast a_k(h_k)dh_1\cdots dh_k
\end{multline*}
Now, the $m>0$-parts of \eqref{classofpairingEZetc} plugs $(Bh)^m(a_0\otimes\cdots\otimes a_{k-2m})$ into this formula, which is of course also well-defined for the broader class of distributions which the application of $Bh$ yields. To see that these terms vanish, we argue that the formula above vanishes for $B$-exact elements Indeed, these $B$-exact elements have the property that they are made up of pure tensors all having one or more $1$'s somewhere (the $1$ here being the Dirac-delta distribution with as value the constant function $1$ on $M$). All the terms where there is a $1$ in the $a_1$ up to $a_k$ slot vanish immediately, since $d(1)=0$. Now for the terms $1\otimes a_1\otimes\cdots\otimes a_k$ we obtain
\begin{multline*}
\frac{(-1)^{\frac{n-k}{2}}}{k!\left(\frac{n-k}{2}\right)!}\int_M{\rm tr}_E(F(\nabla)^{\frac{n-k}{2}})\int_{G^{\times (k-1)}}da_1((h_1\cdots h_{k-1})^{-1})\wedge ((h_1\cdots h_{k-1})^{-1})^\ast da_2(h_1)\wedge\cdots\\
\cdots\wedge (h_{k-1}^{-1})^\ast a_k(h_{k-1})dh_1\cdots dh_{k-1}
\end{multline*}
The integral over $G^k$ becomes one over $G^{k-1}$, since integrating out $a_0((h_1\cdots h_k)^{-1})$ when $a_0=1$, the integral reduces to an integral over the Burgleha space, which is diffeomorphic to $G^{k-1}$. Now, note that since ${\rm tr}_E(F(\nabla)^{\frac{n-k}{2}})$ is closed (it being the differential form inducing the non-equivariant Chern class of $E$), we see that the above formula can be written as
\begin{multline*}
\frac{(-1)^{\frac{n-k}{2}}}{k!\left(\frac{n-k}{2}\right)!}\int_Md\left({\rm tr}_E(F(\nabla)^{\frac{n-k}{2}})\int_{G^{\times (k-1)}}a_1((h_1\cdots h_{k-1})^{-1})\wedge ((h_1\cdots h_{k-1})^{-1})^\ast da_2(h_1)\wedge\cdots\right.\\
\left.\cdots\wedge (h_{k-1}^{-1})^\ast a_k(h_{k-1})dh_1\cdots dh_{k-1}\right)
\end{multline*}
which vanishes since we integrate an exact element.

In conclusion, we see that $c(\Ch_G(E,\nabla))$ only has contributions from the $(m=0)$-parts of \eqref{classofpairingEZetc}, and we conclude:
\begin{multline*}
c(\Ch_G(E,\nabla))(a_0,...,a_k)=\frac{(-1)^{\frac{n-k}{2}}}{k!\left(\frac{n-k}{2}\right)!}\int_M{\rm tr}_E(F(\nabla)^{\frac{n-k}{2}})\int_{G^{\times k}}a_0((h_1\cdots h_k)^{-1})\wedge\\
\wedge((h_1\cdots h_k)^{-1})^\ast da_1(h_1)\wedge\cdots\wedge (h_k^{-1})^\ast a_k(h_k)dh_1\cdots dh_k
\end{multline*}
which equals $\Ch_{\Omega_{E,\nabla}}$.
\end{proof}
\begin{remark}
One would expect the result of \cref{chernsagree} to also hold for non-proper actions. For discrete actions this holds due to work of Gorokhovsky \cite{gorokhovsky}, which we will also discuss below. The main point is to replace $M$ with $M\times X$ where $G$ acts on $X$ properly, and to replace $E\to M$ with $\pr^\ast E\to M\times X$. The essential point in the discrete case is that there exists an algebra map $\Gamma\ltimes \mathsf{A}\to (\Gamma\times\Gamma)\ltimes\mathsf{A}$ when $\Gamma$ is a discrete group and $\mathsf{A}$ is a $\Gamma\times\Gamma$-algebra (with an induced diagonal $\Gamma$-action) which sends $U_gf$ to $U_{g,g}f$. Using this, and K\"unneth-like formulas in periodic cyclic cohomology, one can define an equivariant cup-product, and Gorokhovsky shows that the Chern character of $\pr^\ast E\to M\times X$ is the cup-product of the Chern character of $E\to M$ and the fundamental cycle of $X$.

For the non-discrete case, we expect a parallel argument to work by replacing $M$ by $M\times G/K$ where $K$ is a maximal compact subgroup of $G$ (remark that $G/K$ is a contractible space and the action of $G$ on $G/K$ is proper). We also expect the main problem of defining cup-products can be solved by using results of Nistor \cite{nistor}, which show that localizations of periodic cyclic cohomology of convolution algebras $G\ltimes\mathsf{A}$ are (with a degree shift) isomorphic to the localizations of the periodic cyclic cohomology of $K\ltimes\mathsf{A}$.
\end{remark}

\section{Comparison with known cases}
\subsection{Trivial group actions}
If the group $G$ is the trivial group, the convolution algebra $G\ltimes A$ is simply the commutative algebra $A$. Our double complex in this case becomes
\[\sL(G,A)_{p,q}=(A^+)^{\otimes(p+1)}.\]
However, the vertical differential $d^v\colon \sL(G,A)_{\bullet,q}\to \sL(G,A)_{\bullet,q-1}$ satisfies
\[d^v=\left\{\begin{matrix}
\id &q\text{ is even},\\
0 & q\text{ is odd}.
\end{matrix}\right.\]
In particular, we may replace the total complex $\Tot_\bullet(\sL(G,A))$ by the first row $\sL(G,A)_{\bullet,0}$ and we see that our complex is simply the Hochschild complex of $A^+$.

Similarly, the complex $\sC_{\bullet,\bullet}(G,\Omega_\mathfrak{g}A)$ can be replaced by the complex $\sC_{\bullet,0}(G,\Omega_\mathfrak{g}A)$, which in this case is the `deRham complex of $A$', i.e.\ the mixed complex $(\Omega^\bullet A,0,d)$. The map $\Psi_2\colon\sL(G,A)\to \sC(G,\Omega_\mathfrak{g}A)$ now simply becomes the ordinary HKR map $(A^+)^{\otimes (\bullet+1)}\to \Omega^\bullet A$, and since the map $\EZ\circ\Psi_1\colon \sC_\bullet^\Hoch(G\ltimes A,G\ltimes A)\to \sL_{\bullet,0}(G,A)$ is just the identity, we see that our chain of maps $\Psi$ is simply the HKR-morphism
\[(\sC^\bullet_\Hoch(A,A),b,B)\to(\Omega^\bullet A,0,d).\]

Looking at the case $A=C^\infty_c(M)$ for $M$ a manifold, we note that the equivariant cohomology $\sH_G(M)$ is the deRham cohomology and that Getzler's model $\sC^{\bullet,\bullet}(G,\Omega_\mathfrak{g}C^\infty(M))$ is concentrated in degree $q=0$, where it is given by the deRham-complex of $M$. All in all, the maps
\[c\colon \sH_G^{\text{ev}}(M)\to \sHP^{\dim(M)}(G\ltimes C^\infty_c(M)),\,\, \sH_G^{\text{odd}}(M)\to \sHP^{\dim(M)+1}(G\ltimes C^\infty_c(M))\]
are in this case induced (up to some signs) by the map \[[-]\colon\Omega^\bullet(M)\to\Hom(C^\infty_c(M)^{\dim(m)+1-\bullet},\mathbb{R})\] that takes a differential form $\omega\in\Omega^n(M)$ and sends it so
\[[\omega](f_0,...,f_{\dim(M)-n})=\int_M \omega\wedge f_0df_1\wedge\cdots\wedge df_{\dim(M)-n}.\]
In particular it is the concatenation of the isomorphisms
\[ \sH^\bullet_\dR(M)\xrightarrow{\cong} \sH_{\dim(M)-\bullet}^{\dR,c}(M)\xleftarrow{\cong}\sHP^{\dim(M)-\bullet}(C^\infty_c(M)).\]
Here the first map is an instance of Poincar\'e duality which associates to a $k$-form $\omega\in\Omega^k(M)$, the density $\tilde{\omega}\in \Omega_c^{n-k}(M)^\times$ given by
\[\tilde{\omega}(\alpha)=\int_M\omega\wedge\alpha,\]
and the second map is a consequence of the Hochschild-Kostant-Rosenberg Theorem in the continuous setting.

That the Chern characters are compatible in this case is clear from an immediate determination of the generalized cycle associated to a vector bundle $E\to M$ with a connection $\nabla$. As there is no group action, we see that the underlying curved DGA is the curved DGA
\[\Omega=\Omega(M,\End(E))\]
with differential $d_{\nabla^\End}$ and curvature $F(\nabla)$. The resulting generalized cycle is then seen to be equal (up to a sign) to the current induced by the differential
\[\sum_{i\geq 0}\frac{1}{i!}\tr(F(\nabla)^{\wedge i}),\]
which is exactly the differential form inducing the ordinary Chern character $\Ch(E)\in \sH^\text{ev}_\dR(M)$.
\subsection{Compact groups}
If we put $M=\{{\rm pt}\}$, we recover the convolution algebra of the group $G$ itself. If $G$ is compact, the periodic cyclic cohomology of $C^\infty(G)$ has been computed by Natsume and Nest \cite[\S 1.II]{natsumenest}. In this case, it is concentrated in even degrees, where it is represented by traces $\tau_\phi$ for functions $\phi$ on the spectrum $\hat{G}$ of $G$ which are slowly increasing in a certain sense. Here, for $f\in C^\infty(G)$ this trace $\tau_\phi$ is given by
\begin{equation*}
\tau_\phi(f)=\sum_{\pi\in\widehat{G}}\frac{\phi(\pi)}{\dim(V_\pi)}\int_G f(g)\chi_\pi(g) dg
\end{equation*}
Looking at the equivariant cohomology $\sH^\bullet_G(\cmdtext{pt})$, we obtain from Getzler's model that it is contained in even degrees, where it is given by
\[\sH^{2q}_G(\cmdtext{pt})\cong (\Sym^q(\mathfrak{g}^\ast))^G,\]
the invariant degree $q$ polynomials on $\mathfrak{g}$.

Dissecting the map $c\colon \sH^{\ev}_G(M)\to \sHP^{\ev}(C^\infty(G))$ in this case we notice that the invariant polynomials live in $\sC^{0,0}(G,\Omega_\mathfrak{g}(\{\cmdtext{pt}\}))$, so that the only interesting pairing is with $\sC_{0,0}(G,\Omega_\mathfrak{g}(\{\cmdtext{pt}\}))$. Next, notice that pairing between $\sC^{0,0}(G,\Omega_\mathfrak{g}(\{\cmdtext{pt}\}))$ and $\sC_{0,0}(G,\Omega_\mathfrak{g}(\{\cmdtext{pt}\}))$ kills off polynomials of strictly positive degree as the pairing takes a polynomial $P\in\Sym(\mathfrak{g}^\ast)$ and a function $f\in C^\infty(\mathfrak{g})$ and pairs them by
\[\langle P,f\rangle=P(0)f(0).\]
We conclude the following:
\begin{proposition}The map $c\colon \sH^\ev_G(\cmdtext{pt})\to \sHP^\ev(C^\infty(G))$ takes an invariant polynomial $P\in\Sym(\mathfrak{g}^\ast)\cong \sH^\ev_G(\cmdtext{pt})$ and sends it to the trace $c(P)\in \sHP^0(C^\infty(G))$ given by
\[c(P)(f)=f(e)P(0).\]
\end{proposition}
\begin{remark}
The fact that $\tau(f)=f(e)$ is even a trace on the convolution algebra $C^\infty(G)$ is a consequence of the fact that any compact group is unimodular. Indeed, one checks that
\[\tau([f_1,f_2])=\int_G f_1(g)f_2(g^{-1})dg-\int_G f_2(g)f_1(g^{-1})dg,\]
and these two integrals are equal because of unimodularity.
\end{remark}
\begin{remark}
Under the isomorphism of Natsume and Nest, the trace $c(P)$ corresponds to $\tau_\phi$ for
\[\phi(\pi)=P(0)\dim(V_\pi)^2.\]
This follows from the fact that the character of the regular representation $L^2(G)\cong\widehat{\bigoplus}_{\pi\in\hat{G}}V(\pi)^{\oplus\dim(V_\pi)}$ acts as the Dirac delta distribution at $e\in G$ on the space $C^\infty(G)$.
\end{remark}
Recall that by Proposition \ref{prop-invariantpolygradedtrace} we have a map $\Ch_{\Omega,-}\colon \Sym(\mathfrak{g})^G\to \sHP^\ev(\mathcal{A}_G)$ which sends an invariant polynomial $\gamma$ to the character of the cycle with closed graded trace $\fint_\gamma$. A simple calculation with Fourier inversion shows that under the isomorphism of Natsume and Nest, $\Ch_{\Omega,\gamma}$ is the trace associated to the map
\[\phi(\pi)=\cmdtext{dim}(V_\pi)^2\mathsf{D}_\gamma(\tr(\pi))(0).\]
\subsection{Compact group actions}
In Block-Getzer \cite{BG}, a model for equivariant cyclic homology was presented for when the group $G$ is compact using sheaves over $G$ (with the topology defined by conjugacy-invariant opens) where stalks at $g\in G$ sketch the picture of $M_g=\{p\in M: pg=p\}$ and the Lie algebra $\mathfrak{g}^g$ of the centralizer of $g$ using germs of $G^g$-invariant forms on $M_g$ with polynomial coefficients in $\mathfrak{g}^g$. Using the $C^\infty_{\inv}(G)$-module structure on our complexes, we see that localizing at the identity in our models correspond to the stalk at the identity in the models of Block-Getzler, since we can go from $G$-cohomology to $G$-invariants at no cost by compactness of the group.

As such, we precisely recover the map $\alpha_e$ of \cite[Thm. 3.3]{BG}, from which we infer that
\begin{corollary}\cite[Thm 3.3]{BG}
When the group $G$ is compact, the map $\Psi\colon \sCC(G\ltimes A)\to \Tot(\sCC(G,\Omega_\mathfrak{g}A))$ given by the composition of the diagram \eqref{cdch} is a quasi-isomorphism when localized at the identity.
\end{corollary}
As such, understanding the effect of the map $c\colon \sH^\bullet_G(M)\to \sHP^\bullet(G\ltimes C^\infty_c(M))$ is tantamount to understanding the pairing between $\sC^{\bullet,\bullet}(G,\Omega_\mathfrak{g}C^\infty(M))$ and $\sC_{\bullet,\bullet}(G,\Omega_\mathfrak{g}C_c^\infty(M))$. Since the group $G$ is compact we can the double complexes with the concentration of $G$-cohomology and $G$-homology respectively in their first rows. In particular:
\[ \Tot(\sC^{\bullet,\bullet}(G,\Omega_\mathfrak{g}C^\infty(M))\simeq (\Sym(\mathfrak{g}^\ast)\otimes\Omega^\bullet(M))^G\]
and
\[ \Tot(\sC_{\bullet,\bullet}(G,\Omega_\mathfrak{g}C^\infty_c(M))\simeq (C^\infty(\mathfrak{g})\otimes\Omega_c^\bullet(M))_G.\]

The pairing between these two starts with evaluating at $X=0$ in $\mathfrak{g}$, and the rest is an equivariant instance of the Poincar\'e pairing between differential forms:
\begin{align*}\begin{matrix}
\Omega^{\dim(M)-p}(M)^G\otimes \Omega^p_c(M)_G&\to&\mathbb{R}\\
\langle \omega,\eta\rangle&\mapsto&\int_M \omega\wedge\eta,
\end{matrix}
\end{align*}
which is a (homologically) perfect pairing by Poincar\'e duality. Again we see the story that we almost have a homologically perfect pairing, apart from the fact that we first kill all the behaviour in the $\mathfrak{g}$-direction.
\subsection{Actions of discrete groups}
When the group $G$ is discrete, our constructions and results directly generalize parts of the work done in \cite{gorokhovsky}. In particular, the curved DGA defined in \cref{crvDGA} is precisely the curved DGA defined in \cite[Sect 3]{gorokhovsky} when the group is discrete (of course in this case the moment $\mu$ vanishes).

We also note that when the group is discrete, we obviously overcome the problem where we lose information in the $\mathfrak{g}$-direction when pairing between equivariant cohomology and cyclic homology: indeed, as the group is discrete the Lie algebra $\mathfrak{g}$ is trivial. 

In the discrete case we notice that the convolution algebra $\Gamma\ltimes C^\infty_c(M)$ is a twisted tensor product of the group algebra of $\Gamma$ and the $\Gamma$-algebra $C^\infty_c(M)$. In particular, elements are normally written as sums of elements $U_g f$ for $g\in \Gamma$ and $f\in C^\infty_c(M)$. The product is then given by
\[(U_g f)(U_h f')=U_{gh}f(g\cdot f').\] In his book, Connes \cite[\S III.2.$\delta$]{connes-book} describes a model for the equivariant cohomology $\sH_\Gamma(M)$. In this model, he makes use of maps
\[\gamma\colon  \Gamma^{\times \bullet}\to\Omega_\bullet(M)\]
between products of the group and deRham-currents on $M$. He also gives a map pairing the chains of the $(b,B)$-bicomplex of the convolution algebra $\Gamma\ltimes C^\infty_c(M)$, which -up to signs and combinatorial factors- pairs a map $\gamma$ as above with chains of the convolution algebra by the formula
\[\langle \gamma,(U_{g_0}f_0,...,U_{g_n}f_n)\rangle\sim \gamma(g_0,...,g_n)(f_0 df_1\wedge\cdots\wedge df_n).\]
Furthermore, he shows \cite[Thm III.2.14]{connes-book} that this procedure gives an isomorphism between equivariant cohomology $\sH^\bullet_\Gamma(M)$ and the periodic cyclic cohomology $\sHP^\bullet(\Gamma\ltimes C^\infty_c(M))$ localized at the units.

Using a Poincar\'e duality argument, we can replace currents with differential forms, and we recover Getzler's model for equivariant cohomology. Translating the pairing to this situation and working through the calculations with the Eilenberg-Zilber map, one concludes that our map $c$ between equivariant cohomology $\sH^\bullet_\Gamma(M)$ and periodic cyclic cohomology $\sHP^\bullet(\Gamma\ltimes C^\infty_c(M))$ is precisely the map written down by Connes.

Using this, \cite[Thm 3.1]{gorokhovsky} directly translates to a proof for \cref{chernsagree} in the general case.

\begin{corollary}\cite[Thm 3.1]{gorokhovsky}
When a discrete group $\Gamma$ acts on an oriented manifold $M$, then for any equivariant vector bundle $E$, the two Chern classes $c(\Ch_\Gamma(E)),\Ch_{\Omega_\Gamma(E)}\in \sHP^{\rm even}(\Gamma\ltimes C^\infty_c(M))$ agree.
\end{corollary}
It is noteworthy how the argument simplifies for discrete groups, as opposed to the ideas for a proof we describe in the general case. So let us quickly discuss Gorokhovsky's argument to reduce to the proper case. Again, the main point in reducing to the proper case is replacing $M$ by $M\times X$ such that $\Gamma$ acts properly on $X$.
Next, in \cite[Prop 3.3]{gorokhovsky}, the standard cup product in cyclic cohomology is used to obtain a map
\begin{equation*}
\sHP^{\dim(M)}(\Gamma\ltimes C^\infty_c(M))\otimes \sHP^{\dim(X)}(\Gamma\ltimes C^\infty_c(X))\xrightarrow{\cup} \sHP^{\dim(M\times X)}((\Gamma\times \Gamma)\ltimes C^\infty_c(M\times X))
\end{equation*}
and it is shown that $\Ch_{\Omega(E)}\cup\tau=\Ch_{\Omega(\pr^\ast E)}$. Here $\tau$ is the equivariant orientation class that is associated to the trivial equivariant line bundle over $X$ and $\pr^\ast E$ is seen as a $\Gamma\times \Gamma$-equivariant vector bundle. Then, using a map
\begin{equation*}
\Delta\colon \Gamma\ltimes C^\infty_c(M\times X)\to (\Gamma\times \Gamma)\ltimes C^\infty_C(M\times X)
\end{equation*}
defined by $\Delta(U_gf)=U_{g,g}f$, one obtains a map
\begin{equation*}
\sHP^{\dim(M\times X)}((\Gamma\times \Gamma)\ltimes C^\infty_c(M\times X))\to \sHP^{\dim(M\times X)}(\Gamma\ltimes C^\infty_c(M\times X))
\end{equation*}
and it is shown that this total reduction sends $\Ch_{\Omega(E)}\otimes\tau$ to $\Ch_{\Omega(\pr^\ast E)}$. From that point on, the argument proceeds along roughly the same lines.

Of course, the big difference with the case where $G$ is a non-discrete Lie group is that a map like $\Delta$ does not exist, and there is no obvious way to reduce the cyclic cohomology of $(G\times G)\ltimes A$ (for $A$ an $G\times G$-algebra) to the cyclic cohomology of $G\ltimes A$ where we see $A$ as an $G$-algebra by the diagonal action.

\appendix
\section{Showing that $\Omega_{E,\nabla}$ is an externally curved DGA}
\label{sec:proof}
We use this appendix to prove Proposition \ref{OmegaEcDGA}.
\begin{proposition}
The quadruple $\Omega_{E,\nabla}=(\Omega_E,\ast,D_\nabla,\Theta_\nabla)$ is an externally curved DGA. 
\end{proposition}
\begin{proof}
This proof is essentially the same as the untwisted case, by writing out the explicit equations. The only difference is that we have to navigate the fact that the group action does not commute with the connection, and that $\nabla_{X_M}$ does not equal $\mathcal{L}_X$. For this navigation we have a few claims, most of which follow from a standard type calculation. In what follows $X\in\mathfrak{g}$, $g,h\in G$, $\omega,\eta\in\Omega(M,{\rm End}(E))$ and $\alpha\in C^\infty_c(G,{\rm Sym}\mathfrak{g}\otimes\Omega_c(M,{\rm End}(E)))$.

\textbf{Claim 1:} $\iota_{X_M}g^\ast\omega =g^\ast\iota_{{\rm Ad}_{g^{-1}}(X)_M}\omega$

\textbf{Claim 2:} $\delta(hg)=\delta(h)+h^\ast\delta(g)$

\textbf{Claim 3:} $d_{\nabla^{\rm End}}(g^\ast\omega)=g^\ast d_{\nabla^{\rm End}}\omega+\delta(g)\wedge g^\ast\omega-(-1)^{|\omega|}g^\ast\omega\wedge\delta(g)$

\textbf{Claim 4:} $d^2_{\nabla^{\rm End}}\omega=F(\nabla)\wedge\omega-\omega\wedge F(\nabla)$

\textbf{Claim 5:} $g^\ast F(\nabla)=F(\nabla)-d_{\nabla^{\rm End}}(\delta(g))+\delta(g)\wedge\delta(g)$

\textbf{Claim 6:} $\{d_{\nabla^{\rm End}},\iota_{X_M}\}\omega-\mathcal{L}_X\omega=\mu(X)\wedge\omega-\omega\wedge\mu(X)$

\textbf{Claim 7:} $g^\ast(\mu({\rm Ad}_{g^{-1}}(X)))=\mu(X)-\iota_{X_M}\delta(g)$

\textbf{Claim 8:} $((F(\nabla)+\mu)\delta_e\ast\alpha)(g,X)=(F(\nabla)+\mu(X))\wedge\alpha(g,X)$
\begin{proof}\renewcommand{\qedsymbol}{$\square$ C8}
We use the description of the convolution product as in \ref{conv-distr}. Let $\phi\in C^\infty_c(G,{\rm Sym}\mathfrak{g}^\ast\otimes\Omega_c(M,{\rm End}(E)))$ be a test-function. Fixing the placeholder-variable $g_1$ we need to test the function $g_2\mapsto g_1^{-1}\cdot \phi(g_1g_2)$ on $T_\alpha$, adding the $\mathfrak{g}$-component is this the function sending $(g_2,X)$ to $(g_1^{-1})^\ast\phi(g_1g_2,{\rm Ad}_{g_1}(X))$, so testing against $T_\alpha$ we obtain
\begin{equation*}
\langle T_\alpha(g_2),g_1^{-1}\cdot \phi(g_1g_2)\rangle(X)=\int_G\alpha(g_2,X)\wedge (g_1^{-1})^\ast(\phi(g_1g_2,{\rm Ad}_{g_1}(X)))dg_2
\end{equation*}
and so the function to test against $(F(\nabla)+\mu)\delta_e$ sends $(g_1,X)$ to
\begin{align*}
(g_1\cdot \langle T_\alpha(g_2),g_1^{-1}\cdot \phi(g_1g_2)\rangle)(X)&=g_1^\ast(\langle T_\alpha(g_2),g_1^{-1}\cdot \phi(g_1g_2)\rangle({\rm Ad}_{g_1^{-1}}(X)))\\
&=\int_G (g_1^\ast(\alpha(g_2,{\rm Ad}_{g_1^{-1}}(X)))\wedge \phi(g_1g_2,X)dg_2
\end{align*}
testing this against $(F(\nabla)+\mu)\delta_e$ we plug in $g_1=e$ and take the wedge product to obtain
\begin{equation*}
\langle (F(\nabla)+\mu)\delta_e\ast T_\alpha,\phi\rangle(X)=\int_G (F(\nabla)+\mu(X))\wedge \alpha(g_2,X)\wedge \phi(g_2,X)dg_2
\end{equation*}
and we see that this distribution is precisely of the form $T_\beta$ for
\begin{equation*}
\beta(g,X)=(F(\nabla)+\mu(X))\wedge \alpha(g,X)
\end{equation*}
which proves the claim.
\end{proof}
\textbf{Claim 9:} $(\alpha\ast(F(\nabla)+\mu)\delta_e)(g,X)=\alpha(g,X)\wedge(g^\ast F(\nabla)+g^\ast(\mu({\rm Ad}_{g^{-1}}(X))))$
\begin{proof}\renewcommand{\qedsymbol}{$\square$ C9}
We again let $\phi$ be a test-function and fix the placeholder variable $g_1$. And now test the function $(g_2,X)\mapsto (g_1^{-1})^\ast\phi(g_1g_2,{\rm Ad}_{g_1}(X))$ against $(F(\nabla)+\mu)\delta_e$ obtaining:
\begin{equation*}
\langle (F(\nabla)+\mu)\delta_e(g_2),g_1\cdot \phi(g_1g_2)\rangle(X)=(F(\nabla)+\mu(X))\wedge(g_1^{-1})^\ast\phi(g_1,{\rm Ad}_{g_1}(X))
\end{equation*}
so we obtain
\begin{equation*}
(g_1\cdot \langle (F(\nabla)+\mu)\delta_e(g_2),g_1\cdot \phi(g_1g_2)\rangle)(X)=(g_1^\ast F(\nabla)+g_1^\ast(\mu({\rm Ad}_{g_1^{-1}}(X))))\wedge \phi(g_1,X)
\end{equation*}
testing this against $T_\alpha(g_1)$ we obtain
\begin{equation*}
\langle T_\alpha\ast (F(\nabla)+\mu)\delta_e,\phi\rangle(X)=\int_G \alpha(g_1,X)\wedge (g_1^\ast F(\nabla)+g_1^\ast(\mu({\rm Ad}_{g_1^{-1}}(X))))\wedge\phi(g_1,X)dg_1
\end{equation*}
and this is precisely of the form $T_\beta$ for
\begin{equation*}
\beta(g,X)=\alpha(g,X)\wedge(g^\ast F(\nabla)+g^\ast(\mu({\rm Ad}_{g^{-1}}(X))))
\end{equation*}
which proves the claim
\end{proof}
\textbf{Claim 10:} $d_{\nabla^{\rm End}}(\omega\wedge\eta)=d_{\nabla^{\rm End}}(\omega)\wedge\eta+(-1)^{|\omega|}\omega\wedge(d_{\nabla^{\rm End}}\eta)$.

\textbf{Claim 11:} $\left.\frac{d}{dt}\right|_{t=0}\delta(e^{tX})=-d_{\nabla^{\rm End}}(\mu(X))-\iota_{X_M}F(\nabla)$

Using these claims we can show that $\Omega_{E,\nabla}$ is an externally curved DGA. First note that the multiplication is associative by general machinery since it is the convolution product induced by a $G$-algebra.

We then show that $D_\nabla(\alpha\ast\beta)=(D_\nabla(\alpha))\ast\beta)+(-1)^{|\alpha|}\alpha\ast(D_\nabla(\beta))$ for all $\alpha,\beta\in\Omega$. We investigate the three terms of which $D_\nabla$ is made up:
\begin{align*}
d_{\nabla^{\rm End}}((\alpha\ast\beta)(g,X))&=\int_G d_{\nabla^{\rm End}}(\alpha(h,X)\wedge h^\ast\beta(h^{-1}g,{\rm Ad}_{h^{-1}}(X)))dh.
\end{align*}
Using Claim 10 this divides up into
\begin{align*}
d_{\nabla^{\rm End}}((\alpha\ast\beta)(g,X))=&\int_G (d_{\nabla^{\rm End}}(\alpha(h,X)))\wedge h^\ast\beta(h^{-1}g,{\rm Ad}_{h^{-1}}(X))dh\\
&+(-1)^{|\alpha|}\int_G\alpha(h,X)\wedge (d_{\nabla^{\rm End}}(h^\ast\beta(h^{-1}g,{\rm Ad}_{h^{-1}}(X))))dh.
\end{align*}
Using Claim 3 on the second line results in
\begin{align*}
d_{\nabla^{\rm End}}((\alpha\ast\beta)(g,X))=&\int_G (d_{\nabla^{\rm End}}(\alpha(h,X)))\wedge h^\ast\beta(h^{-1}g,{\rm Ad}_{h^{-1}}(X))dh\\
&+(-1)^{|\alpha|}\int_G\alpha(h,X)\wedge h^\ast(d_{\nabla^{\rm End}}(\beta(h^{-1}g,{\rm Ad}_{h^{-1}}(X)))dh\\
&+(-1)^{|\alpha|}\int_G\alpha(h,X)\wedge\delta(h)\wedge h^\ast\beta(h^{-1}g,{\rm Ad}_{h^{-1}}(X))dh\\
&+(-1)^{|\alpha|+|\beta|}\int_G\alpha(h,X)\wedge h^\ast\beta(h^{-1}g,{\rm Ad}_{h^{-1}}(X))\wedge\delta(h)dh.
\end{align*}
And then using Claim 2 on the $\delta(h)$-term on the last line to replace it by $h^\ast\delta(h^{-1}g)-\delta(g)$ and we obtain
\begin{align*}
d_{\nabla^{\rm End}}((\alpha\ast\beta)(g,X))=&\int_G (d_{\nabla^{\rm End}}(\alpha(h,X)))\wedge h^\ast\beta(h^{-1}g,{\rm Ad}_{h^{-1}}(X))dh\\
&+(-1)^{|\alpha|}\int_G\alpha(h,X)\wedge h^\ast(d_{\nabla^{\rm End}}(\beta(h^{-1}g,{\rm Ad}_{h^{-1}}(X)))dh\\
&+(-1)^{|\alpha|}\int_G\alpha(h,X)\wedge\delta(h)\wedge h^\ast\beta(h^{-1}g,{\rm Ad}_{h^{-1}}(X))dh\\
&+(-1)^{|\alpha|+|\beta|}\int_G\alpha(h,X)\wedge h^\ast(\beta(h^{-1}g,{\rm Ad}_{h^{-1}}(X))\wedge\delta(h^{-1}g))dh\\
&+(-1)^{|\alpha|+|\beta|+1}\int_G\alpha(h,X)\wedge h^\ast\beta(h^{-1}g,{\rm Ad}_{h^{-1}}(X))\wedge\delta(g)dh.
\end{align*}
Note that the last line here cancels to the $\wedge\delta$-part of $D_\nabla$ when applied to $\alpha\ast\beta$.

Then working on the $\iota_{X_M}$-term we have
\begin{align*}
\iota_{X_M}((\alpha\ast\beta)(g,X))=&\int_G\iota_{X_M}(\alpha(h,X)\wedge h^\ast\beta(h^{-1}g,{\rm Ad}_{h^{-1}}(X)))dh\\
=&\int_G(\iota_{X_M}\alpha(h,X))\wedge h^\ast\beta(h^{-1}g,{\rm Ad}_{h^{-1}}(X)))dh\\
&+(-1)^{|\alpha|}\int_G\alpha(h,X)\wedge \iota_{X_M}(h^\ast\beta(h^{-1}g,{\rm Ad}_{h^{-1}}(X)))dh.
\end{align*}
Applying Claim 1 in the second line we obtain
\begin{align*}
\iota_{X_M}((\alpha\ast\beta)(g,X))=&\int_G(\iota_{X_M}\alpha(h,X))\wedge h^\ast\beta(h^{-1}g,{\rm Ad}_{h^{-1}}(X)))dh\\
&+(-1)^{|\alpha|}\int_G\alpha(h,X)\wedge h^\ast(\iota_{{\rm Ad}_{h^{-1}}(X)_M}\beta(h^{-1}g,{\rm Ad}_{h^{-1}}(X)))dh.
\end{align*}
Combining all this, we see
\footnotesize
\begin{align*}
D_\nabla((\alpha\ast\beta))(g,X)=&\int_G\left(d_{\nabla^{\rm End}}(\alpha(h,X))+(-1)^{|\alpha|}\alpha(h,X)\wedge\delta(h)+\iota_{X_M}\alpha(h,X)\right)\wedge h^\ast\beta(h^{-1}g,{\rm Ad}_{h^{-1}}(X))dh\\
&+(-1)^{|\alpha|}\int_G\alpha(h,X)\wedge h^\ast\left(d_{\nabla^{\rm End}}(\beta(h^{-1}g,{\rm Ad}_{h^{-1}}(X))+\right.\\
&\left.(-1)^{|\beta|}\beta(h^{-1}g,{\rm Ad}_{h^{-1}}(X))\wedge\delta(h^{-1}g)+\iota_{{\rm Ad}_{h^{-1}}(X)_M}\beta(h^{-1}g,{\rm Ad}_{h^{-1}}(X))\right)dh\\
=&\int_G (D_\nabla\alpha)(h,X)\wedge h^{-1}\beta(h^{-1}g,{\rm Ad}_{h^{-1}}(X))dh\\
&+(-1)^{|\alpha|}\int_G \alpha(h,X)\wedge h^{-1}((D_\nabla\beta)(h^{-1}g,{\rm Ad}_{h^{-1}}(X))dh\\
=&\,\,((D_\nabla\alpha)\ast\beta)(g,X)+(-1)^{|\alpha|}(\alpha\ast(D_\nabla\beta))(g,X).
\end{align*}
\normalsize
Lastly we show that $D_\nabla^2=[\Theta_\nabla,-]$. First we calculate $D_\nabla^2$:
\begin{align*}
(D_\nabla^2\alpha)(g,X)=&\,\,d_{\nabla^{\rm End}}((D_\nabla\alpha)(g,X))-(-1)^{|\alpha|}(D_\nabla\alpha)(g,X)\wedge\delta(g)+\iota_{X_M}((D_\nabla\alpha)(g,X))\\
=&\,\,d_{\nabla^{\rm End}}^2(\alpha(g,X))+(-1)^{|\alpha|}d_{\nabla^{\rm End}}(\alpha(g,X)\wedge \delta(g))+d_{\nabla^{\rm End}}(\iota_{X_M}(\alpha(g,X)))\\
&-(-1)^{|\alpha|}d_{\nabla^{\rm End}}(\alpha(g,X))\wedge\delta(g)-\alpha(g,X)\wedge\delta(g)\wedge\delta(g)\\
&-(-1)^{|\alpha|}\iota_{X_M}(\alpha(g,X))\wedge\delta(g)+\iota_{X_M}(d_{\nabla^{\rm End}}(\alpha(g,X)))\\
&+(-1)^{|\alpha|}\iota_{X_M}(\alpha(g,X)\wedge\delta(g))+\iota_{X_M}\iota_{X_M}\alpha(g,X).
\end{align*}
Now we use Claim 10, the fact that $\iota_{X_M}$ is a graded derivation and $\iota_{X_M}^2=0$ to obtain
\begin{align*}
(D_\nabla^2\alpha)(g,X)=&\,\,d^2_{\nabla^{\rm End}}(\alpha(g,X))+\alpha(g,X)\wedge d_{\nabla^{\rm End}}(\delta(g))-\alpha(g,X)\wedge\delta(g)\wedge\delta(g)\\
&+\{d_{\nabla^{\rm End}},\iota_{X_M}\}(\alpha(g,X))+\alpha(g,X)\wedge\iota_{X_M}(\delta(g)).
\end{align*}
Now we use Claim 4 and Claim 6 to obtain
\begin{align*}
(D^2_\nabla\alpha)(g,X)=&(F(\nabla)+\mu(X))\wedge\alpha(g,X)+\mathcal{L}_X(\alpha(g,X))\\
&-\alpha(g,x)\wedge\left(F(\nabla)-d_{\nabla^{\rm End}}(\delta(g))+\delta(g)\wedge\delta(g)+\mu(X)-\iota_{X_M}\delta(g)\right).
\end{align*}
Then Claims 5 and 7 give us
\begin{multline*}
(D^2_\nabla\alpha)(g,X)=(F(\nabla)+\mu(X))\wedge\alpha(g,X)+\mathcal{L}_X(\alpha(g,X))\\
-\alpha(g,X)\wedge(g^\ast F(\nabla)+g^\ast\mu({\rm Ad}_{g^{-1}}(X))).
\end{multline*}
To finish, note that the calculation in the untwisted case carries over to here to give $\mathcal{L}_X=[\Theta,-]$, which combined with Claims 8 and 9 will give
\begin{equation*}
(D^2_\nabla\alpha)(g,X)=[\Theta_\nabla,\alpha](g,X).
\end{equation*}

Lastly we check the Bianchi identity $D_\nabla(\Theta_\nabla\ast\alpha)=\Theta_\nabla\ast D_\nabla\alpha$. For notational clarity we give the three terms of $D_\nabla$ explicit names, namely $\nabla$, $\delta$ and $\iota$ for
\[
(\nabla\alpha)(g,X)=d_{\nabla^{\rm End}}(\alpha(g,X)),
\qquad
(\delta\alpha)(g,X)=(-1)^{|\alpha|}\alpha(g,X)\wedge\delta(g,X),\]
\[
(\iota\alpha)(g,X)=\iota_{X_M}\alpha(g,X).
\]

We calculate
\begin{align*}
\nabla(\Theta_\nabla\ast\alpha)(g,X)=&\,\,d_{\nabla^{\rm End}}(\mathcal{L}_X(\alpha(g,X)))-\left.\frac{d}{dt}\right|_{t=0}d_{\nabla^{\rm End}}(\alpha(e^{tX}g,X))\\
&+d_{\nabla^{\rm End}}(F(\nabla)\wedge \alpha(g,X))+d_{\nabla^{\rm End}}(\mu(X)\wedge\alpha(g,X))
\end{align*}
\begin{align*}
(\Theta_\nabla\ast\nabla\alpha)(g,X)=&\,\,\mathcal{L}_X(d_{\nabla^{\rm End}}(\alpha(g,X)))-\left.\frac{d}{dt}\right|_{t=0} d_{\nabla^{\rm End}}(\alpha(e^{tX}g,X))\\
&+F(\nabla)\wedge d_{\nabla^{\rm End}}(\alpha(g,X))+\mu(X)\wedge d_{\nabla^{\rm End}}
\end{align*}
On the nose the second terms cancel and combining that $d_{\nabla^{\rm End}}$ is a graded derivation for the wedge-product and the classical Bianchi-identity $d_{\nabla^{\rm End}}(F(\nabla))=0$ so do the third terms. We obtain
\begin{equation*}
(\nabla(\Theta_\nabla\ast\alpha)-\Theta_\nabla\ast\nabla\alpha)(g,X)=[d_{\nabla^{\rm End}},\mathcal{L}_X](\alpha(g,X))+(d_{\nabla^{\rm End}}\mu(X))\wedge\alpha(g,X)
\end{equation*}
Next we calculate the terms that appear here, first we have
\begin{align*}
\mathcal{L}_X(d_{\nabla^{\rm End}}(\alpha(g,X)))=&\,\,\left.\frac{d}{dt}\right|_{t=0} e^{tX}\cdot d_{\nabla^{\rm End}}(\alpha(g,X))\\
({\rm C3})\to\,\,\,=&\,\,\left.\frac{d}{dt}\right|_{t=0}d_{\nabla^{\rm End}}(e^{tX}\cdot\alpha(g,X))-\left.\frac{d}{dt}\right|_{t=0}\delta(e^{tX})\wedge e^{tX}\cdot \alpha(g,X)\\
&+(-1)^{|\alpha|}\left.\frac{d}{dt}\right|_{t=0}(e^{tX}\cdot\alpha(g,X))\wedge \delta(e^{tX})
\end{align*}
Using the fact that $\delta(e)=0$ this results in
\begin{equation*}
[d_{\nabla^{\rm End}},\mathcal{L}_X](\alpha(g,X))=\left(\left.\frac{d}{dt}\right|_{t=0}\delta(e^{tX})\right)\wedge\alpha(g,X)-(-1)^{|\alpha|}\alpha(g,X)\wedge\left(\left.\frac{d}{dt}\right|_{t=0}\delta(e^{tX})\right)
\end{equation*}
Using Claim 11 it now follows that
\begin{equation*}
(\nabla(\Theta_\nabla\ast\alpha)-\Theta_\nabla\ast\nabla\alpha)(g,X)=-\iota_{X_M}F(\nabla)\wedge \alpha(g,X)-(-1)^{|\alpha|}\alpha(g,X)\wedge \left(\left.\frac{d}{dt}\right|_{t=0}\delta(e^{tX})\right)
\end{equation*}
Next we look at the $\delta$-part and see
\begin{align*}
\delta(\Theta_\nabla\ast\alpha)(g,X)=&(-1)^{|\alpha|}\mathcal{L}_X(\alpha(g,X))\wedge\delta(g)-(-1)^{|\alpha|}\left.\frac{d}{dt}\right|_{t=0}\alpha(e^{tX}g,X)\wedge\delta(g)\\
&+(-1)^{|\alpha|}F(\nabla)\wedge\alpha(g,X)\wedge\delta(g)+(-1)^{|\alpha|}\mu(X)\wedge\alpha(g,X)\wedge\delta(g)
\end{align*}
\begin{align*}
(\Theta_\nabla\ast\delta\alpha)(g,X)=&(-1)^{|\alpha|}\mathcal{L}_X(\alpha(g,X)\wedge\delta(g))-(-1)^{|\alpha|}\left.\frac{d}{dt}\right|_{t=0}\alpha(e^{tX}g,X)\wedge\delta(e^{tX}g)\\
&+(-1)^{|\alpha|}F(\nabla)\wedge\alpha(g,X)\wedge\delta(g)+(-1)^{|\alpha|}\mu(X)\wedge\alpha(g,X)\wedge\delta(g)
\end{align*}
From which we conclude
\begin{align*}
(\delta(\Theta_\nabla\ast\alpha)-(\Theta_\nabla\ast\delta\alpha))(g,X)&=(-1)^{|\alpha|}\alpha(g,X)\wedge\left(\left.\frac{d}{dt}\right|_{t=0}(\delta(e^{tX}g)-e^{tX}\cdot\delta(g))\right)\\
&=(-1)^{|\alpha|}\alpha(g,X)\wedge\left(\left.\frac{d}{dt}\right|_{t=0}\delta(e^{tX})\right)
\end{align*}
Lastly we take a look at the $\iota$-part.
\begin{align*}
\iota(\Theta_\nabla\ast\alpha)(g,X)=&\iota_{X_M}(\mathcal{L}_X(\alpha(g,X)))-\left.\frac{d}{dt}\right|_{t=0}\iota_{X_M}(\alpha(e^{tX}g,X))\\
&+\iota_{X_M}(F(\nabla)\wedge\alpha(g,X))+\iota_{X_M}(\mu(X)\wedge\alpha(g,X))
\end{align*}
\begin{align*}
(\Theta_\nabla\ast\iota\alpha)(g,X)=&\mathcal{L}_X(\iota_{X_M}(\alpha(g,X)))-\left.\frac{d}{dt}\right|_{t=0}\iota_{X_M}(\alpha(e^{tX}g,X))\\
&+F(\nabla)\wedge\iota_{X_M}(\alpha(g,X))+\mu(X)\wedge\iota_{X_M}(\alpha(g,X))
\end{align*}
Here we find
\begin{equation*}
(\iota(\Theta_\nabla\ast\alpha)-\Theta_\nabla\ast\iota\alpha)(g,X)=[\iota_{X_M},\mathcal{L}_X](\alpha(g,X))+\iota_{X_M}F(\nabla)\wedge\alpha(g,X)
\end{equation*}
Exactly similar to the untwisted case we have $[\iota_{X_M},\mathcal{L}_X]=0$, and so we see that the contributions from $\nabla$, $\delta$ and $\iota$ cancel and we conclude that the Bianchi identity
\begin{equation*}
D_\nabla(\Theta_\nabla\ast\alpha)=\Theta_\nabla\ast D_\nabla\alpha
\end{equation*}
indeed holds.

We conclude that $\Omega_{E,\nabla}$ is an externally curved DGA, as claimed.
\end{proof}

\bibliographystyle{alpha}
\bibliography{references}

\end{document}